\input ./style/arxiv-general.cfg
\documentclass[aap,MSNbibl,dvips]{arximspdf}
\makeatletter
   \@ifpackageloaded{graphicx}{}{\usepackage{graphicx}}
\makeatother
\doi{10.1214/14-AAP1079}% Updated by VTEXPTS2LaTeX.exe, 22.12.2014
\volume{25}
\issue{6}
\pubyear{2015}
\firstpage{3434}
\lastpage{3464}
\docsubty{FLA}

\makeatletter
\newcommand{\rrvert}{\vert}
\newcommand{\rrVert}{\Vert}
\newcommand{\llvert}{\vert}
\newcommand{\llVert}{\Vert}
\newtheorem{theorem}{Theorem}
\newtheorem{lemma}{Lemma}[section]
\newproclaim{rem}{Remark}
\renewcommand{\mid}{|}
\def\ZZ{\mathbb{Z}}
\def\var{\operatorname{var}}
\makeatother

\begin{document}
\begin{frontmatter}

%\dochead{}
\title{Coexistence of grass, saplings and trees in the Staver--Levin
forest model}
\runtitle{Coexistence in the Staver--Levin forest model}

\begin{aug}
% Corresponding author: Rick Durrett - rtd@math.duke.edu% Updated by
%VTEXPTS2LaTeX.exe, 22.12.2014 15:08
\author[A]{\fnms{Rick}~\snm{Durrett}\corref{}\thanksref{T1}\ead[label=e1]{rtd@math.duke.edu}}
\and
\author[A]{\fnms{Yuan}~\snm{Zhang}\thanksref{T1}\ead[label=e2]{yzhang@math.duke.edu}}
\runauthor{R. Durrett and Y. Zhang}
\affiliation{Duke University}
%\dedicated{}
\address[A]{Department of Mathematics\\
Duke University\\
Box 90320 \\
Durham, North Carolina 27708-0320\\
USA\\
\printead{e1}\\
\phantom{E-mail:\ }\printead*{e2}}
\end{aug}
\thankstext{T1}{Supported in part by NSF Grant DMS-10-05470 from the
probability program.}

% HISTORY:
%
\received{\smonth{8} \syear{2014}}% Updated by VTEXPTS2LaTeX.exe,
%22.12.2014 15:08
%\revised{\smonth{} \syear{}}

% ABSTRACT
%
\begin{abstract}
In this paper, we consider two attractive stochastic spatial models in
which each site can be in state 0, 1 or 2: Krone's model in which
0${}={}$vacant, 1${}={}$juvenile and 2${}={}$a mature individual
capable of giving birth, and the Staver--Levin forest model in which
0${}={}$grass,
1${}={}$sapling and 2${}={}$tree. Our first result shows that if
$(0,0)$ is an unstable fixed point of the mean-field ODE for
densities of 1's and 2's then when the range of interaction is large,
there is positive probability of survival starting from a finite set
and a stationary distribution
in which all three types are present. The result we obtain in this way
is asymptotically sharp for Krone's model.
However, in the Staver--Levin forest model, if $(0,0)$ is attracting
then there may also be another stable fixed point for the ODE,
and in some of these cases there is a nontrivial stationary distribution.
\end{abstract}

% KEYWORDS
% Pirmas kwd is didziosios raides
%
\begin{keyword}[class=AMS]
\kwd[Primary ]{60K35}
\kwd{82C22}
\kwd[; secondary ]{82B43}
\end{keyword}
\begin{keyword}
\kwd{Coexistence}
\kwd{stationary distribution}
\kwd{percolation}
\kwd{block construction}
\end{keyword}
\end{frontmatter}

%s1 #&#
\section{Introduction}\label{sec1}

In a recent paper published in Science \cite{SR0}, Carla Staver, Sally
Archibald and Simon Levin argued that tree cover does not increase
continuously with rainfall but rather is constrained to low ($<$50\%,
``savanna'') or high ($>$75\%, ``forest'') levels. In follow-up work
published in Ecology \cite{SR1}, the American Naturalist \cite{SR2}
and Journal of Mathematical Biology \cite{SSL}, they studied the
following ODE for the evolution of the fraction of land covered by
grass $G$, saplings $S$ and trees $T$:
%
%
%e1 #&#
\begin{eqnarray}
 \frac{dG}{dt}&=&\mu S+\nu T- \beta GT,
\nonumber
\\
 \frac{dS}{dt}&=&\beta GT-\omega(G) S-\mu S, \label{mODE}
\\
 \frac{dT}{dt}&=&\omega(G) S-\nu T.
\nonumber
\end{eqnarray}
Here, $\mu\ge\nu$ are the death rates for saplings and trees, and
$\omega(G)$ is the rate at which saplings grow into trees. Fires
decrease this rate of progression, and the incidence of fires
is an increasing function of the fraction of grass, so $\omega(G)$ is
decreasing.
Studies suggest (see \cite{SR2} for references) that regions with
tree cover below about 40\% burn frequently but fire is rare above this
threshold, so they used
an $\omega$ that is close to a step function.

The ODE in (\ref{mODE}) has very interesting behavior: it may have two
stable fixed points, changing the values of parameters may lead to Hopf
bifurcations, and if the system has an extra type of savanna trees,
there can be periodic orbits.
In this paper, we will begin the study of the corresponding spatial
model. The state at time $t$ is $\chi_t\dvtx \ZZ^d \to\{0,1,2\}$, where
$0 ={}$grass, $1 ={}$sapling and $2 ={}$tree. Given the application, it would be
natural to restrict our attention to $d=2$, but since the techniques we
develop will be applicable to other systems we consider the general case.

In the forest model, it is natural to assume that dispersal of seeds is
long range.
To simplify our calculations, we will not use a continuous dispersal
distribution for tree seeds, but instead let $f_i(x,L)$ denote the
fraction of sites of type $i$ in the box $x + [-L,L]^d$ and declare
that site $x$ changes:
\begin{itemize}
\item$0\rightarrow1$ at rate $\beta f_2(x,L)$,
\item$1\rightarrow2$ at rate $\omega(f_0(x,\kappa L))$,
\item$1\rightarrow0$ at rate $\mu$,
\item$2\rightarrow0$ at rate $\nu$.
\end{itemize}

The configuration with all sites 0 is an absorbing state. This
naturally raises the question of finding conditions that guarantee that
coexistence occurs, i.e., there is a stationary distribution in
which all three types are present. Our model has three states but it is
``attractive,'' that is, if $\chi_0(x) \le\chi'_0(x)$ for all $x$
then we can construct the two processes on the same space so that this
inequality holds for all time. From this, it follows from the usual
argument that if we start from $\chi^2_0(x) \equiv2$ then $\chi^2_t$
converges to a limit $\chi^2_\infty$ that is a translation invariant
stationary distribution, and there will be a nontrivial stationary
distribution if and only if $P( \chi^2_\infty(0) = 0 ) < 1$. Since
2's give birth to 1's and 1's grow into 2's, if $\chi^2_\infty$ is
nontrivial then both species will be present with positive density in
$\chi^2_\infty$.

If $\omega\equiv\gamma$ is constant, $\mu=1+\delta$, and $\nu=1$,
then our system reduces to one studied by Krone \cite{KR}.
In his model, 1's are juveniles who are not yet able to reproduce.
Krone proved the existence of nontrivial stationary distributions in
his model by using a simple comparison between
the sites in state 2 and a discrete time finite-dependent oriented
percolation. In the percolation process, we have an edge from $(x,n)
\to(x+1,n+1)$ if a 2 at $x$ at time $n\varepsilon$ will give birth to a
1 at $x+1$, which then grows to a 2 before time $(n+1)\varepsilon$,
and there are no deaths at $x$ or $x+1$ in $[n\varepsilon
,(n+1)\varepsilon
]$. As the reader can imagine, this argument produces a very crude
result about
the parameter values for which coexistence occurs.

A simple comparison shows that if we replace the decreasing function
$\omega(G)$ in the Staver--Levin model, $\chi_t$, by the constant
$\omega=\omega(1)$, to obtain
a special case $\eta_t$ of Krone's model, then $\chi_t$ dominates
$\eta_t$ in the sense that given $\chi_0 \ge\eta_0$ the two processes
can be coupled so that $\chi_t \ge\eta_t$
for all $t$. Because of this, we can prove existence of nontrivial
stationary distribution in the Staver--Levin model by
studying Krone's model. To do this under the assumption of long range
interactions, we begin with the mean field ODE:
%
%
%e2 #&#
\begin{eqnarray}
\label{KmODE}  \frac{dG}{dt}&=&\mu S+\nu T- \beta GT,\nonumber
\\
 \frac{dS}{dt}&=&\beta GT-(\omega+\mu) S,
\\
 \frac{dT}{dt}&=&\omega S-\nu T.\nonumber
\end{eqnarray}
Here, $\omega$ is a constant. When it is a function, we will write
$\omega(G)$.

Since $G+S+T=1$, we can set $G=1-S-T$ and reduce the system to two
equations for $S$ and $T$.
To guess a sufficient condition for coexistence in the long range
limit we note that:
%

%le1.1 #&#
\begin{lemma}
For the mean-field ODE (\ref{KmODE}), $S=0,T=0$ is not an attracting
fixed point if
%
%
%e3 #&#
\begin{equation}
\label{ConSur} \mu\nu<\omega(\beta-\nu).
\end{equation}
On the other hand, $S=0,T=0$ is attracting if
%
%
%e4 #&#
\begin{equation}
\label{ConExt} \mu\nu>\omega(\beta-\nu).
\end{equation}
\end{lemma}

\begin{pf} When $(S,T) \approx(0,0)$, and hence $G \approx1$, the
mean-field ODE is approximately
\[
\pmatrix{ dS/dt
\vspace*{3pt}\cr
dT/dt }\approx A \pmatrix{ S
\vspace*{3pt}\cr
T } \qquad\mbox{where } A =
\pmatrix{ -(\omega+\mu) & \beta
\vspace*{3pt}\cr
\omega& -\nu}. %
\]
The trace of $A$, which is the sum of its eigenvalues is negative, so
$(0,0)$ is not attracting if the determinant of $A$,
which is the product of the eigenvalues is negative. Since $(\omega
+\mu)\nu-\beta\omega<0$ if and only if $\mu\nu<(\beta-\nu
)\omega$,
we have proved the desired result. Similarly, $(0,0)$ is attracting if
the determinant of $A$ is positive, which implies (\ref{ConExt}).
\end{pf}

%
%th1 #&#
\begin{theorem}\label{mainth}\label{teo1}
Let $\eta_t$ be Krone's model with parameters that satisfy (\ref{ConSur}).
Then when $L$ is large enough, $\eta_t$ survives with positive
probability starting from a finite number of
nonzero sites and $\eta_t$ has a nontrivial stationary distribution.
\end{theorem}

Foxall \cite{Fox} has shown that for Krone's model the existence of
nontrivial stationary distribution is equivalent to survival for a
finite set of nonzero sites, so we only have to prove one of these
conclusions. However, our proof is via a block construction, so we get
both conclusions at the same time.

Our next result is a converse, which does not require the assumption of
long range.

%
%th2 #&#
\begin{theorem} \label{Kdie}\label{teo2}
Suppose $\mu\nu\ge\omega(\beta-\nu)$. Then for any $L>0$, Krone's
model $\eta_t$ dies out, that is, for any initial configuration $\eta_0$
with finitely many nonzero sites:
\[
\lim_{t\rightarrow\infty}P\bigl(\eta_t(x)\equiv0\bigr)=1.
\]
If $\mu\nu> \omega(\beta-\nu)$, then for any initial
configuration $\eta_0$ and any $x\in\ZZ^d$, the probability
\[
P\bigl(\exists t>0, \mbox{ s.t. } \eta_s(x)=0, \forall s\ge t\bigr)=1.
\]
\end{theorem}

The second conclusion implies that there is no nontrivial
stationary distribution.
Comparing with Krone's model, we see that if $\mu\nu> \omega
(0)(\beta-\nu)$ then the Staver--Levin model dies out.

%s1.1 #&#
\subsection{Survival when zero is stable}\label{sec11}
When $\mu\nu>\omega(1)(\beta-\nu)$, the Staver--Levin ODE (\ref
{mODE}) may have another stable fixed point in the positive density
region (and also an unstable fixed point in between), the Staver--Levin
model, like the quadratic contact process studied by \cite{DNeu,Neu}
and \cite{BesD} may have a nontrivial stationary
distribution when $(0,0)$ is attracting.

Based on the observation in \cite{SR1} mentioned above, it is natural
to assume that $\omega(\cdot)$ is a step function. In our proof, we let
%
%
%e5 #&#
\begin{equation}
\label{RoG} \omega(G) = \cases{ \omega_0, &\quad$G\in[0,1-
\delta_0)$,
\vspace*{3pt}\cr
\omega_1, &\quad$G\in[1-
\delta_0,1]$,}
\end{equation}
where $\omega_0>\omega_1$ and $\delta_0\in(0,1)$. However,
according to the monotonicity of $\chi_t$, our result about the
existence of a nontrivial stationary distribution will also hold if one
replaces the ``$=$'' in $(\ref{RoG})$ by ``$\ge$'' since the new
process dominates the old one.

To prove the existence of a nontrivial stationary distribution under
the assumption of long range, a natural approach would be to show that
when $L\rightarrow\infty$
and space is rescaled by diving by $L$, the Staver--Levin model
converges weakly to the solution of following integro-differential equation:
%
%
%e6 #&#
\begin{eqnarray}\label{IDE}
\frac{dS(x,t)}{dt}&=&\beta D^T_1(x,t) G-\mu S-\omega
\bigl[D^G_\kappa(x,t)\bigr]S,
\nonumber\\[-8pt]\\[-8pt]\nonumber
\frac{dT(x,t)}{dt}&=&\omega\bigl[D^G_\kappa(x,t) \bigr]S-\nu T,
\end{eqnarray}
where $G = 1-S-T$ and
\begin{eqnarray*}
 D^T_1(x,t) &=&\frac{\int_{x+[-1,1]^d}T(y,t) \,dy_1\cdots dy_d}{2^d},
\\
 D^G_\kappa(x,t)&=&\frac{\int_{x+[-\kappa,\kappa]^d}G(y,t)
\,dy_1\cdots dy_d}{(2\kappa)^d},
\end{eqnarray*}
are
the local densities of trees and grass on the rescaled lattice. The
first problem with this approach is that since the density is computed
by examining all sites in a square, there is not a good dual process,
which was the key to proofs in \cite{SG,DNeu,Neu}
and \cite{BesD}. The second problem is that one does not know much
about the limiting IDE. Results of Weinberger \cite{WE} show the
existence of wave speeds and provide a convergence theorem in the case
of a single equation, but we do not know of results for a pair of equations.

To avoid these difficulties, we will construct test functions
$S_{\mathrm{test}}$ and $T_{\mathrm{test}}$, so that under (\ref{IDE}), the derivatives
will always be positive for all $x$ in $\overline{\{T_{\mathrm{test}}>0\}}$ and
$\overline{\{S_{\mathrm{test}}>0\}}$, where $\bar A$ stands for the closure of
set $A$. The positive derivative implies that after a positive time the
solution will dominate translates of the initial condition by positive
and negative amounts. Monotonicity then implies that the solution will
expand linearly and the result follows from a block construction.
Details can be found in Section~\ref{sec8}.

%
%th3 #&#
\begin{theorem}\label{SurHD}\label{teo3}
Recall the definition of $\omega(G)$ in (\ref{RoG}). Under condition
%
%
%e7 #&#
\begin{equation}
\label{SFsur2} \beta\omega_0>2^{d}\nu(\mu+
\omega_0),
\end{equation}
there is a constant $S_0(\beta,\omega_0,\nu,\mu)$ so
that if
\[
\delta_0 \in\bigl(0,2^{-d}S_0\bigr) %
\]
the Staver--Levin forest model $\chi_t $ survives when $L$ is large
for any $\omega_1 \ge0$ and $\kappa> 0$.
\end{theorem}

Combining Theorems \ref{teo1}--\ref{teo3}, we have the following results for the
Staver--Levin model:
\begin{longlist}[2.]
\item[1.] When $\mu\nu<\omega(1)(\beta-\nu)$, $\chi_t$ survives from
a finite set of nonzero sites when $L$ is large.

\item[2.] When $\mu\nu\ge\omega(0)(\beta-\nu)$, $\chi_t$ dies out
from a finite set of nonzero sites for all $L \ge1$.

\item[3.] When $\mu\nu\ge\omega(1)(\beta-\nu)$, under the hypotheses
of Theorem \ref{SurHD}, $\chi_t$ can still survive from a finite set
of nonzero sites when $L$ is large, no matter how small is~$\omega_1$.
\end{longlist}

%s1.2 #&#
\subsection{Sketch of the proof of Theorem \texorpdfstring{\protect\ref{mainth}}{1}}\label{sec12}

Most of the remainder of the paper is devoted to the proof of Theorem~\ref{teo1}. We will now describe the main ideas and then explain where the
details can be found.

\begin{longlist}[(iii)]
\item[(i)] \textit{The key idea} is due to Grannan and Swindle \cite
{GS}. They consider a model of a catalytic surface in which atoms of
type $i = 1,2$
land at vacant sites (0's) at rate $p_i$, while adjacent $1,2$ pairs
turn into $0,0$ at rate $\infty$. If after a landing event,
several $1,2$ pairs are created, one is chosen at random to be removed.
The first type of event is the absorption
of an atom onto the surface of the catalyst, while the second is a
chemical reaction, for example, carbon monoxide $\mathrm{CO}$
and oxygen~$\mathrm{O}$ reacting to produce $\mathrm{CO}_2$. The last reaction occurs in
the catalytic converted in your car,
but the appropriate model for that system is more complicated. An
oxygen molecule $\mathrm{O}_2$ lands and dissociates to two $\mathrm{O}$
bound to the surface when a pair of adjacent sites is vacant. See
Durrett and Swindle \cite{DS94} for more details
about the phase transition in the system.

Suppose without loss of generality that
$p_1+p_2=1$. In this case, Grannan and Swindle \cite{GS} showed that
if $p_1 \neq p_2$ the only possible stationary distributions concentrate
configurations that are $\equiv1$ or $\equiv2$. Mountford and Sudbury
\cite{MS} later improved this result
by showing that if $p_1 > 1/2$ and the initial configuration has
infinitely many 1's then the system converges
to the all 1's state.

The key to the Grannan--Swindle argument was to consider
\[
Q(\eta_t) = \sum_x e^{-\lambda\llVert x\rrVert } q
\bigl[\eta_t(x)\bigr], %
\]
where $\llVert x\rrVert = \sup_i \llvert x_i\rrvert $ is the
$L^\infty$ norm, $q(0)=0$,
$q(1)=1$, and $q(2)=-1$.
If $\lambda$ is small enough then $dEQ/dt \ge0$ so $Q$ is a bounded
submartingale,
and hence converges almost surely to a limit. Since an absorption or
chemical reaction in $[-K,K]^d$
changes $Q$ by an amount $\ge\delta_K$, it follows that such events
eventually do not occur.

\item[(ii)] \textit{Recovery from small density} is the next step. We will pick
$\varepsilon_0>0$ small, let $\ell= [\varepsilon_0L]$ be the integer
part of $\varepsilon_0L$
and divide space into \textit{small boxes} $\hat{B}_x = 2\ell x + (-\ell
, \ell]^d$. To make the number of 1's and 2's in the various small boxes
sufficient to describe the state of the process, we declare two small
boxes to be neighbors if all of their
points are within an $L^\infty$ distance $L$. For the ``truncated
process,'' which is stochastically bounded by $\eta_t$, and in which
births of trees can only occur between
sites in neighboring small boxes, we will show that if $\kappa\in
(d/2,d)$ and we start with
a configuration that has $L^\kappa$ nonzero sites in $\hat{B}_0$ and
0 elsewhere, then the system will recover and produce
a small box $\hat{B}_x$ at time $\tau$ in which the density of
nonzero sites is $a_0>0$ and $P( \tau> t_0 \log L ) < L^{d/2-\kappa}$.
See Lemma \ref{reclem}.
To prove this, we use an analogue of Grannan and Swindle's $Q$. The
fact that $(0,0)$ is an unstable fixed point implies
$dEQ/dt > 0$ as long as the density in all small boxes is ${\le} a_0$.

\item[(iii)] \textit{Bounding the location of the positive density box} is the
next step.
To do this, we use a comparison with branching random walk to show that
the small box $\hat{B}_x$ with density $a_0$
constructed in step~(ii) is not too far from 0. Random walk estimates will
later be used to control how far it will wander as we iterate the
construction. For this step, it is
important that the truncated process is invariant under reflection, so
the mean displacement is 0. If we try to work directly with the original
interacting particle system $\eta_t$ then it is hard to show that the
increments between box locations are independent and have mean 0. It is
for this reason we introduced the truncated process.

\item[(iv)] \textit{Moving particles}. The final ingredient in the block
construction is to show that given a small block
$\hat B_x$ with positive density and any $y$ with $\llVert y-x
\rrVert _1 \le[c
\log L]$ then if c is small enough it is very likely that there will\vspace*{1pt}
be ${\ge} L^\kappa$ particles in $\hat B_y$ at time $[c \log L]$.
Choosing $y$ appropriately and then using the recovery
lemma, we can get lower bounds on the spread of the process.

\item[(v)] \textit{Block construction}. Once we have completed steps (ii), (iii)
and (iv), it is straightforward to show that our system dominates a
one-dependent oriented percolation. This shows that the system survives
from a finite set with positive probability and proves the existence of
nontrivial stationary distribution.
\end{longlist}

The truncated process is defined in Section~\ref{sec2} and a graphical
representation is used to couple it, Krone's model and the
Staver--Levin model. In Section~\ref{sec3}, we use the Grannan--Swindle argument
to do step (ii). The dying out result, Theorem~\ref{teo2}, is proved in
Section~\ref{sec4}. In Sections~\ref{sec5}, \ref{sec6}, and \ref{sec7}, we take care of steps (iii), (iv)
and (v). In Section~\ref{sec8}, we prove Theorem~\ref{teo3}.

%s2 #&#
\section{Box process and graphical representation}\label{sec2}

For some fixed $\varepsilon_0>0$ which will be specified in (\ref
{ep0def}), let $l=[\varepsilon_0 L]$ and divide space $\ZZ^d$ into
small boxes:
\[
\hat{B}_x=2lx+(-l,l]^d, \qquad x\in\ZZ^d.
\]
For any $x\in\ZZ^d$, there is a unique $x'$ such that $x\in\hat
{B}_{x'}$. Define the new neighborhood of interaction as follows: for
any $y\in\hat{B}_{y'}$,
$y\in\mathcal{N}(x)$ if and only if
\[
\sup_{z_1\in\hat{B}_{x'}, z_2\in\hat{B}_{y'}}\llVert z_1-z_2\rrVert\le
L. %
\]
It is easy to see that $\mathcal{N}(x)\subset B_x(L)$ where $B_x(L)$
is the $L^\infty$ neighborhood centered at $x$ with range $L$. To show that
%
%e8 #&#
\begin{equation}
\mathcal{N}(x) \supset B_x\bigl((1-4\varepsilon_0)L\bigr)
\label{NxLB}
\end{equation}
we note that if $\llVert x-y\rrVert \le(1-4\varepsilon_0)L$, $z_1\in
\hat{B}_{x'}$,
and $z_2\in\hat{B}_{y'}$, where $\hat{B}_{x'}$ and $\hat{B}_{y'}$
are the small boxes containing $x$ and $y$:
\[
\llVert z_1-z_2\rrVert\le\llVert z_1-x
\rrVert+\llVert x-y\rrVert+\llVert y-z_2\rrVert\le4\bigl([
\varepsilon_0 L]\bigr)+(1-4\varepsilon_0)L\le L. %
\]
Given the new neighborhood $\mathcal{N}(x)$, we define the truncated
version of Krone's model $\xi_t$ by its transition rates:
\begin{center}\vspace*{6pt}
\begin{tabular}{cc}
transtion & at rate \\
$1\rightarrow0$ & $\mu$ \\
$2\rightarrow0$ & $\nu$ \\
$1\rightarrow2$ & $\omega$ \\
$0\rightarrow1$ & $\beta N_2 (\mathcal{N}(x) )/(2L+1)^d$,
\end{tabular}
\end{center}\vspace*{6pt}
where $N_i(S)$ stands for the number of $i$'s in the set $S$.

For any $x\in\ZZ^d$ and $\xi\in\{0,1,2\}^{\ZZ^d}$, define
$n_i(x,\xi)$ to be the number of type $i$'s in the small box $\hat
{B}_{x}$ in the configuration $\xi$.
The box process is defined by
\[
\zeta_t(x)= \bigl(n_1(x,\xi_t),n_2(x,
\xi_t) \bigr)\qquad \forall x\in\ZZ^d.\vadjust{\goodbreak}
\]
Then $\zeta_t$ is a Markov process on $\{(n_1,n_2)\dvtx n_1,n_2\ge0,
n_1+n_2\le\llvert \hat{B}_0\rrvert \}^{\ZZ^d}$ in which

\begin{center}\vspace*{6pt}
\begin{tabular}{lc}
\hphantom{$\zeta_t(x)$} transition & at rate \\
$\zeta_t(x)\rightarrow\zeta_t(x) - (1,0)$ & $\mu\zeta^1_t(x)$ \\[2pt]
$\zeta_t(x)\rightarrow\zeta_t(x) - (0,1)$ & $\nu\zeta^2_t(x)$ \\[2pt]
$\zeta_t(x)\rightarrow\zeta_t(x) + (-1,1)$ & $\omega\zeta^1_t(x)$\\[2pt]
$\zeta_t(x)\rightarrow\zeta_t(x) + (1,0)$ & $\zeta^0_t(x)\sum_{y\dvtx\hat{B}_y\subset\mathcal{N}(x)} \beta\zeta^2_t(y)$,
\end{tabular}
\end{center}\vspace*{6pt}
where
\[
\zeta^0_t(x)=\llvert\hat{B}_0\rrvert-
\zeta^1_t(x)-\zeta^2_t(x)
\]
be the number of 0's in that small box.

Because $\zeta_t$ only records the number of particles in any small
box, and the neighborhood
is defined so that all sites in the same small box have the same
neighbors, the distribution of $\zeta_t$ is symmetric
under reflection in any axis. The main use for this observation is that
the displacement of the location of the positive density
box produced by the recovery lemma in Section~\ref{sec3} has mean 0.

%s2.1 #&#
\subsection{Graphical representation}\label{sec2.1}

We will use the graphical representation similar as in \cite{KR} to
construct Krone's model $\eta_t$ and the truncated version $\xi_t$ on
the same probability space, so that:

\begin{longlist}[($\ast\ast$)]
\item[($\ast$)] If $\eta_0 \ge\xi_0$, then we will have $\eta_t\ge\xi
_t$ for all $t$.
\end{longlist}

Note that $\mu\ge\nu$. We use independent families of Poisson
processes for each $x\in\ZZ^d$, as follows:
\begin{enumerate}
\item[$\{V_n^x\dvtx n\ge1\}$] with rate $\nu$. We put an $\times$ at space--time
point $(x,V_n^x)$ and write a $\delta_{12}$ next to it to indicate a
death will occur if $x$ is occupied by a 1 or a 2.

\item[$\{U_n^x\dvtx n\ge1\}$] with rate $\mu-\nu$. We put an $\times$ at
space--time point $(x,U_n^x)$ and write a $\delta_{1}$ next to it to
indicate a death will occur if $x$ is occupied by 1.

\item[$\{W_n^x\dvtx n\ge1\}$] with rate $\omega$. We put an $\bullet$ at
space--time point $(x,W_n^x)$ which indicates that if $x$ is in state 1,
it will become a 2.

\item[$\{T^{x,y}_n\dvtx n\ge1 \}$] with rate $\beta/\llvert B_0\rrvert $
for all $y\in
\mathcal{N}(x)$. We draw a solid arrow from $(x,T_n^{x,y})$ to
$(y,T_n^{x,y})$ to indicate that if $x$ is occupied by a 2 and $y$ is
vacant, then a birth will occur at $x$ in either process.

\item[$\{T^{x,y}_n\dvtx n\ge1\}$] with rate $\beta/\llvert B_0\rrvert $ for
all $y\in
B_x(L)-\mathcal{N}(x)$. We draw a dashed arrow from $(x,T_n^{x,y})$ to
$(y,T_n^{x,y})$ to indicate that if $x$ is occupied by a 2 and $y$ is
vacant then a birth will occur at $x$ in the process $\xi_t$.
\end{enumerate}
Standard arguments that go back to Harris \cite{HT} over
forty years ago guarantee that we have constructed the desired
processes. Since each flip preserves $\eta_s \ge\xi_s$, the
stochastic order ($\ast$) is satisfied.

To finish the construction of the Staver--Levin model, $\chi_t$, we
add another family of Poisson process $\{ \hat W_n^x\dvtx n \ge1\}$ with
rate $1-\omega$, and independent random variables $w_{x,n}$ uniform on
$(0,1)$. At any time $\hat W_n^x$ if $x$ is in state 1, it will
increase to state 2 if
\[
w_{n,x} > \frac{ \omega(f_0(x,\kappa L)) - \omega}{1-\omega}.
\]
These events take care of the extra growth of 1's into 2's in $\chi
_t$. Again every flip preserves $\chi_s \ge\eta_s$ so we have:

\begin{longlist}
\item[($\ast\ast$)] If $\chi_0 \ge\eta_0$ then we will have $\chi_t\ge
\eta_t$ for all $t$.
\end{longlist}

%s3 #&#
\section{Recovery lemma}\label{sec3}

Given (\ref{ConSur}), one can pick a $\theta$, which must be ${>}1$,
such that
%
%
%e9 #&#
\begin{equation}
\label{Theta} \frac{\mu+\omega}{\omega}<\theta<\frac{\beta}{\nu}
\end{equation}
so we have
\[
\theta\omega-(\omega+\mu)>0, \qquad\beta-\theta\nu>0 %
\]
and since the inequalities above are strict, we can pick some $a_0>0$
and $\rho\in(0,1)$ such that
%
%
%e10 #&#
\begin{equation}
\label{PDrift} \theta\omega-(\omega+\mu)\ge\rho, \qquad\beta(1-4a_0)-
\theta\nu\ge\theta\rho.
\end{equation}
Now we can let the undetermined $\varepsilon_0$ in the definition of
$\xi
_t$ in Section~\ref{sec2} be a positive constant such that
%
%e11 #&#
\begin{equation}
\label{ep0def} (1-4\varepsilon_0)^d>1-2a_0.
\end{equation}
Fix some $\alpha\in(d/2,d)$. We start with an initial configuration
in $\Xi_0$, the $\xi_0$ that have $\xi_0(x)=0$ for all $x \notin
\hat B_0$
and the number of nonzero sites in $\hat B_0$ is at least $L^\alpha$.
We define a stopping time $\tau$:
%
%
%e12 #&#
\begin{equation}
\label{ToR} \tau=\inf\bigl\{t\dvtx \exists x\in\ZZ^d\mbox{ such that }
n_1(x,\xi_t)+n_2(x,\xi_t)\ge
a_0\llvert\hat{B}_0\rrvert\bigr\}.
\end{equation}

%
%le3.1 #&#
\begin{lemma}[(Recovery lemma)] \label{reclem}
Suppose we start the truncated version of Krone's
model from a $\xi_0 \in\Xi_0$. Let $t_0=2d/\rho$. When $L$ is large,
%
%
%e13 #&#
\begin{equation}
\label{PoS1} P(\tau>t_0\log L)<L^{d/2-\alpha}.
\end{equation}
\end{lemma}

\begin{pf}
As mentioned in the \hyperref[sec1]{Introduction}, we consider
\[
Q(\xi_t)=\lambda^d\sum_{x\in\ZZ^d}
e^{-\lambda\llVert x\rrVert }w\bigl[\xi_t(x)\bigr],
\]
where $\lambda=L^{-1}a_0/2$ and
\[
w\bigl[\xi(x)\bigr]=\cases{ 0, &\quad if $\xi(x)=0$,
\vspace*{3pt}\cr
1, &\quad if $
\xi(x)=1$,
\vspace*{3pt}\cr
\theta, &\quad if $\xi(x)=2$.}
\]
If we imagine $\mathbb{R}^d$ divided into cubes with centers at
$\lambda\ZZ^d$
and think about sums approximating an integral, then we see that
%
%e14 #&#
\begin{equation}
\lambda^d\sum_{x\in\ZZ^d} e^{-\lambda\llVert x\rrVert }
\le e^{\lambda/2} \int_{\mathbb{R}^d} e^{-\llVert z\rrVert } \,dz \le
U_{\fontsize{8.36pt}{11pt}{(\ref{UBsum})}} \label{UBsum}
\end{equation}
for all $\lambda\in(0,1]$. From this, it follows that
%
%e15 #&#
\begin{equation}
Q(\xi_t) \le\theta U_{\fontsize{8.36pt}{11pt}{(\ref{UBsum})}}. \label{UBsum1}
\end{equation}

%re1 #&#
\begin{rem}
Here, and in what follows, we subscript important
constants by the lemmas or formulas where they were first introduced,
so it will be easier for the reader to find where they are defined.
$U$'s are upper bounds that are independent of $\lambda\in(0,1]$.
\end{rem}

Our next step toward Lemma \ref{reclem} is to study the infinitesimal mean
\[
\mu(\xi)=\lim_{\delta t\downarrow0}\frac{E[Q(\xi_{t+\delta
t})-Q(\xi_t)\mid\xi_t=\xi]}{\delta t}. %
\]

%
%le3.2 #&#
\begin{lemma} \label{infmean}
For all $\xi$ such that $n_1(x,\xi)+n_2(x,\xi)\le a_0\llvert \hat
{B}_0\rrvert $
and for all \mbox{$x\in\ZZ^d$},
$\mu(\xi)\ge\rho Q(\xi)$ where $\rho$ is defined in (\ref{PDrift}).
\end{lemma}

\begin{pf} Straightforward calculation gives
%
%e16 #&#
\begin{eqnarray}
\label{muf} \frac{\mu(\xi)}{\lambda^d} &=&\sum_{\xi(x)=1}
\bigl[(\theta-1)\omega-\mu\bigr]e^{-\lambda\llVert x\rrVert }
\nonumber\\[-8pt]\\[-8pt]\nonumber
&&{}+\sum_{\xi(x)=0}\beta\frac{N_2[\mathcal
{N}(x)]}{(2L+1)^d}e^{-\lambda\llVert x\rrVert }
-\sum_{\xi(x)=2}\theta\nu e^{-\lambda\llVert x\rrVert }.
\end{eqnarray}
For the second term in the equation above, we interchange the roles of
$x$ and $y$ then rearrange the sum:
\[
\sum_{\xi(x)=0} \beta\frac{N_2[\mathcal
{N}(x)]}{(2L+1)^d}e^{-\lambda\llVert x\rrVert }=
\sum_{\xi(x)=2}(2L+1)^{-d}\sum
_{y\in\mathcal{N}(x), \xi(y)=0} \beta e^{-\lambda\llVert y\rrVert }. %
\]
Noting that $\lambda=L^{-1}a_0/2$, and that for any $x$ and $y\in
\mathcal{N}(x)\subset B_x$, $-L\le\llVert y\rrVert -\llVert
x\rrVert \le L$, we have
\[
e^{-\lambda\llVert y\rrVert }\ge e^{-a_0/2}e^{-\lambda\llVert
x\rrVert }\ge(1-a_0)e^{-\lambda\llVert x\rrVert }.
\]
Using this with $n_1(x,\xi)+n_2(x,\xi)<a_0\llvert \hat{B}_0\rrvert $, and
$B_x[(1-4\varepsilon_0)L]\subset\mathcal{N}(x)$ from~(\ref{NxLB}),
\begin{eqnarray*}
\sum_{\xi(x)=0} \beta\frac{N_2[\mathcal
{N}(x)]}{(2L+1)^d}e^{-\lambda\llVert x\rrVert }&
\ge& (1-a_0)\sum_{\xi(x)=2}\beta
\frac{N_0(B_x[(1-4\varepsilon_0)L])}{(2L+1)^d}e^{-\lambda\llVert
x\rrVert }
\\
&\ge& (1-a_0)\bigl[(1-4\varepsilon_0)^d-a_0
\bigr]\sum_{\xi(x)=2}\beta e^{-\lambda
\llVert x\rrVert }.
\end{eqnarray*}
Recall that by (\ref{ep0def}), $\varepsilon_0$ is small enough so that
$(1-4\varepsilon_0)^d>1-2a_0$. This choice implies
\begin{eqnarray*}
\sum_{\xi(x)=0} \beta\frac{N_2[\mathcal
{N}(x)]}{(2L+1)^d}e^{-\lambda\llVert x\rrVert }&>&(1-a_0)
(1-3a_0)\sum_{\xi
(x)=2}\beta e^{-\lambda\llVert x\rrVert }
\\
&>& (1-4a_0)\sum_{\xi(x)=2}\beta
e^{-\lambda\llVert x\rrVert }.
\end{eqnarray*}

Combining inequality above with (\ref{muf}) and (\ref{PDrift}) gives
\begin{eqnarray*}
\mu(\xi)&\ge&\lambda^d \sum_{\xi(x)=1}\bigl[(
\theta-1)\omega-\mu\bigr]e^{-\lambda\llVert x\rrVert } +\lambda^d
\sum
_{\xi(x)=2}\bigl[(1-4a_0)\beta-\theta\nu
\bigr]e^{-\lambda\llVert
x\rrVert }
\\
&\ge&\rho Q(\xi),
\end{eqnarray*}
which proves the desired result.
\end{pf}

Then for any initial configuration $\xi_0$, define
%
%
%e17 #&#
\begin{equation}
\label{DyM} M_t=Q(\xi_t)-Q(\xi_0)-\int
_0^t\mu(\xi_s)\,ds.
\end{equation}
According to Dynkin's formula, $M_t$ is a martingale with $EM_t=0$.

%
%le3.3 #&#
\begin{lemma} \label{M2tbd}\label{le3.3}
There are constants $L_{\fontsize{8.36pt}{11pt}{\ref{M2tbd}}}$ and $U_{\fontsize{8.36pt}{11pt}{\ref{M2tbd}}}<\infty$
so that when $L\ge L_{\fontsize{8.36pt}{11pt}{\ref{M2tbd}}}$, we have $EM_t^2 \le U_{\fontsize{8.36pt}{11pt}{\ref{M2tbd}}} L^{-d}t$ for all $t\ge0$, and hence
%
%
%e18 #&#
\begin{equation}
\label{VoMax} E \Bigl(\sup_{s\le t}M_s^2
\Bigr)\le4U_{\fontsize{8.36pt}{11pt}{\ref{M2tbd}}} L^{-d} t.
\end{equation}
\end{lemma}

\begin{pf}
Using (\ref{muf}) and (\ref{UBsum}), we see that
%
%e19 #&#
\begin{equation}
\bigl\llvert\mu(\xi_t)\bigr\rrvert\le C^{(1)}_{\fontsize{8.36pt}{11pt}{\ref{M2tbd}}}=
\theta(\beta+\omega+\mu+\nu)U_{\fontsize{8.36pt}{11pt}{(\ref{UBsum})}}. \label{UBmu}
\end{equation}
To calculate $EM_t^2$, let $t^n_i = it/n$.
%
%e20 #&#
\begin{eqnarray}\label{EM2t}
EM_t^2 &=& \sum_{i=0}^{n-1}E(M_{t^n_{i+1}}-M_{t^n_i})^2
\nonumber\\[-8pt]\\[-8pt]\nonumber
&=& \sum_{i=0}^{n-1} E \biggl[Q(
\xi_{t^n_{i+1}})-Q(\xi_{t^n_i}) -\int_{t^n_i}^{t^n_{i+1}}
\mu(\xi_s)\,ds \biggr]^2.
\end{eqnarray}

The path of $M_s, s\in[0,t]$ is always a right continuous function
with left limit. To control the limit of the sum, we first consider the
total variation of $M_s$, $s\in[0,t]$. For each $n$, let
\[
V^{(n)}_t=\sum_{i=0}^{n-1}
\llvert M_{t^n_{i+1}}-M_{\xi_{t^n_i}}\rrvert. %
\]
By definition,
\begin{eqnarray*}
V^{(n)}_t&\le&\sum_{i=0}^{n-1}
\bigl\llvert Q(\xi_{t^n_{i+1}})-Q(\xi_{t^n_i})\bigr\rrvert+\sum
_{i=0}^{n-1} \biggl\llvert\int
_{t^n_i}^{t^n_{i+1}} \mu(\xi_s)\,ds\biggr\rrvert
\\
&\le& V_t+\int_0^t \bigl\llvert
\mu(\xi_s)\bigr\rrvert ds\le V_t+ C^{(1)}_{\fontsize{8.36pt}{11pt}{\ref{M2tbd}}}
t,
\end{eqnarray*}
where
\[
V_t=\sum_{s\in\Pi_t}\bigl\llvert Q(
\xi_s)-Q(\xi_{s-})\bigr\rrvert%
\]
to be the total variation of $Q(\xi_s)$ in $[0,t]$ and $\Pi_t$ be the
set of jump times of $\xi_s$ in $[0,t]$, which is by definition a
countable set. To control $V_t$, write $\Pi_t=\bigcup_{k=0}^\infty\Pi
^{(k)}_t$, where for each $k$, $\Pi^{(k)}_t$ is the set of times in
which $\xi$ has a transition at a vertex contained in
$H_k=B_0(kL)\setminus B_0((k-1)L)$. Then according to (\ref{UBsum}),
there is some $L_{\fontsize{8.36pt}{11pt}{\ref{M2tbd}}}<\infty$ and $C^{(2)}_{\fontsize{8.36pt}{11pt}{\ref{M2tbd}}},
C^{(3)}_{\fontsize{8.36pt}{11pt}{\ref{M2tbd}}}<\infty$, such that for all $L\ge L_{\fontsize{8.36pt}{11pt}{\ref{M2tbd}}}$,
%
%e21 #&#
\begin{eqnarray}
\label{UBTV} %
EV_t&\le&\sum_{k=0}^\infty
\Bigl[ E \bigl(\bigl\llvert\Pi^{(k)}_t\bigr\rrvert\bigr)
\cdot\sup_{x\in H_k}\mathop{\sup_{
\xi,\xi'\in\{0,1,2\}^{Z^d}:}}_{\xi(y)=\xi'(y)\ \forall y\neq x}
\bigl\llvert Q(\xi_s)-Q(\xi_{s-})\bigr\rrvert\Bigr]
\nonumber\\[-8pt]\\[-8pt]\nonumber
&\le& C^{(2)}_{\fontsize{8.36pt}{11pt}{\ref{M2tbd}}}L^{-d} \sum
_{k=0}^\infty k^{d-1} e^{-\lambda k}t\le
C^{(3)}_{\fontsize{8.36pt}{11pt}{\ref{M2tbd}}}t
\end{eqnarray}
which implies that $V_t<\infty$ almost surely, and that $M_t$ is a
process with finite variation and definitely bounded. Using Proposition
3.4 on page 67 of \cite{EKZ} and the fact that $M_t$ is a bounded
right-continuous martingale, we have
%
%e22 #&#
\begin{equation}
\label{CtQv} \sum_{i=0}^{n-1}E(M_{t^n_{i+1}}-M_{t^n_i})^2
\stackrel{L^1} {\longrightarrow} [M]_t,
\end{equation}
where $[M]_t$ is the quadratic variation of $M_t$. Noting that for any $n$
\[
\sum_{i=0}^{n-1}E(M_{t^n_{i+1}}-M_{t^n_i})^2
\equiv EM_t^2, %
\]
combining this with the $L^1$ convergence in (\ref{CtQv}), we have
\[
EM_t^2=E[M]_t. %
\]
Since $M_t$ is a martingale of finite variation, Exercise 3.8.12 of
\cite{B02} implies
%
%e23 #&#
\begin{equation}
\label{QvJ} [M]_t=\sum_{s\in\Pi_t}
\bigl(Q(\xi_s)-Q(\xi_{s-}) \bigr)^2.
\end{equation}
So for $E[M]_t$, similar as in (\ref{UBTV}), there is some $U_{\fontsize{8.36pt}{11pt}{\ref{M2tbd}}}<\infty$, such that when $L\ge L_{\fontsize{8.36pt}{11pt}{\ref{M2tbd}}}$
%
%e24 #&#
\begin{eqnarray}
\label{UBQV} E[M]_t&\le&\sum_{k=0}^\infty
\Bigl[ E \bigl(\bigl\llvert\Pi^{(k)}_t\bigr\rrvert\bigr)
\cdot\sup_{x\in H_k}\mathop{\sup_{
\xi,\xi'\in\{0,1,2\}^{Z^d}:}}_{ \xi(y)=\xi'(y)\ \forall y\neq x
}\bigl(Q(\xi_s)-Q(\xi_{s-}) \bigr)^2 \Bigr]
\nonumber\\[-8pt]\\[-8pt]\nonumber
&\le& C^{(2)}_{\fontsize{8.36pt}{11pt}{\ref{M2tbd}}}L^{-2d} \sum
_{k=0}^\infty k^{d-1} e^{-2\lambda k}t\le
U_{\fontsize{8.36pt}{11pt}{\ref{M2tbd}}}L^{-d}t.
\end{eqnarray}
Equation (\ref{UBQV}) immediately implies that
\[
EM_t^2=E[M]_t\le U_{\fontsize{8.36pt}{11pt}{\ref{M2tbd}}}L^{-d}t
\]
which completes the proof.
\end{pf}

At this point, we have all the tools needed in the proof of Lemma \ref{reclem}.
If $\xi_0 \in\Xi_0$, there is a $u_{\fontsize{8.36pt}{11pt}{\ref{reclem}}}>0$ such that for
all $\xi_0$ in Lemma \ref{reclem}:
\[
u_{\fontsize{8.36pt}{11pt}{\ref{reclem}}}L^{-d+\alpha}\le Q(\xi_0). %
\]
Using (\ref{VoMax}) now
\[
E \Bigl(\sup_{s\le t_0\log L}M_s^2 \Bigr)
\le4U_{\fontsize{8.36pt}{11pt}{\ref{M2tbd}}} L^{-d} t_0 \log L %
\]
so by Chebyshev's inequality and the fact that $\alpha>d/2$:
\begin{eqnarray*}
P \Bigl( \sup_{s\le t_0\log L}\llvert M_s\rrvert\ge
u_{\fontsize{8.36pt}{11pt}{\ref{reclem}}}L^{-d+\alpha
}/2 \Bigr) &\le&\frac{8 U_{\fontsize{8.36pt}{11pt}{\ref{M2tbd}}}L^{-d} t_0\log L}{u_{\fontsize{8.36pt}{11pt}{\ref{reclem}}}^2L^{-2(d-\alpha)}}
\\
&=&O\bigl(L^{-2\alpha+d}\log L\bigr)
\\
&=&o\bigl(L^{d/2-\alpha}\bigr) \to0.
\end{eqnarray*}

Consider the event $\{\tau>t_0\log L\}$. For any $s\le t_0\log L$,
$n_1(x,\xi_s)+\break n_2(x,\xi_s)<a_0\llvert \hat{B}_0\rrvert $,
for all $x\in\ZZ^d$, so by Lemma \ref{infmean}, $\mu(\xi_s)\ge
\rho Q(\xi_s)$. Consider the set
\[
A = \Bigl\{ \sup_{s\le t_0\log L}\llvert M_s\rrvert<
u_{\fontsize{8.36pt}{11pt}{\ref{reclem}}}L^{-d+\alpha}/2 \Bigr\} \cap\{\tau>t_0\log L\}.
\]
On $A$, we will have that for all $t\in[0,t_0\log L]$,
\[
Q(\xi_t)\ge u_{\fontsize{8.36pt}{11pt}{\ref{reclem}}}L^{-d+\alpha}/2 + \rho\int
_0^t Q(\xi_s)\,ds. %
\]
If we let $f(t)=e^{\rho t}u_{\fontsize{8.36pt}{11pt}{\ref{reclem}}}L^{-d+\alpha}/2$, then
\[
f(t)= u_{\fontsize{8.36pt}{11pt}{\ref{reclem}}}L^{-d+\alpha}/2+\rho\int_0^t
f(s)\,ds. %
\]
Reasoning as in the proof of Gronwall's inequality:

%
%le3.4 #&#
\begin{lemma}
On the event $A$, $Q(\xi_t)\ge f(t)$ for all $t\in[0,t_0\log L]$.
\end{lemma}

\begin{pf} Suppose the lemma does not hold. Let $t_1=\inf\{t\in
[0,t_0\log L]\dvtx\break  Q(\xi_t)<f(t)\}$. By right-continuity of $Q(\xi_t)$,
$Q(\xi_{t_1})\le f(t_1)$ and $t_1>0$. However, by definition of $t_1$,
we have $Q(\xi_{t})\ge f(t)$ on $[0,t_1)$, and by right-continuity of
$Q(\xi_t)$ near $t=0$, the inequality is strict in a neighborhood of
0. Thus, we have
\begin{eqnarray*}
Q(\xi_{t_1})&\ge&\frac{u_{\fontsize{8.36pt}{11pt}{\ref{reclem}}}L^{-d+\alpha}}{2}+\rho\int_0^{t_1}
Q(\xi_s)\,ds
\\
&>&\frac{u_{\fontsize{8.36pt}{11pt}{\ref{reclem}}}L^{-d+\alpha}}{2}+\rho\int_0^{t_1}
f(s)\,ds=f(t_1)
\end{eqnarray*}
which is a contradiction to the definition of $t_1$.
\end{pf}

Recalling that $t_0=2d/\rho$
\[
f(t_0\log L)= e^{\rho t_0\log L}u_{\fontsize{8.36pt}{11pt}{\ref{reclem}}}L^{-d+\alpha
}/2=u_{\fontsize{8.36pt}{11pt}{\ref{reclem}}}L^{d+\alpha}/2.
\]
When $L$ is large, this will be $\ge\theta U_{\fontsize{8.36pt}{11pt}{(\ref{UBsum})}}$, the
largest possible value of $Q(\xi_t)$. Thus, the assumption that
$P(A>0)$ has lead to a contradiction,
and we have completed the proof of Lemma \ref{reclem}.
\end{pf}

%s4 #&#
\section{Proof of Theorem \texorpdfstring{\protect\ref{Kdie}}{2}}\label{sec4}

As in the proof of Lemma \ref{infmean}, we are able to prove the
extinction result in Theorem \ref{Kdie}, which does not require the
assumption of long range.

\begin{pf} %{Proof of Theorem \protect\ref{Kdie}}
When $\mu\nu\ge\omega(\beta-\nu)$, if $\beta\le\nu
$, the system dies out since $\eta_t$ can be bounded by a subcritical
contact process with birth rate $\beta$ and death rate $\nu$ (the
special case of $\eta_t$ when $\omega=\infty$). Otherwise, we can
find a $\theta'$ such that
\[
\frac{\mu+\omega}{\omega}\ge\theta'\ge\frac{\beta}{\nu}>1. %
\]
For $\eta_t$ starting from $\eta_0$ with a finite number of nonzero
sites, consider
\[
S(\eta_t)=\sum_{x\in\ZZ^d}1_{\eta_t(x)=1}+
\theta'1_{\eta_t(x)=2}. %
\]
Similarly, let $\mu(\eta_t)$ be the infinitesimal mean of $S(\eta
_t)$. Repeating the calculation in the proof of Lemma \ref{infmean},
we have
\[
\mu(\eta_t)=\sum_{x\in\ZZ^d}\bigl[\omega
\bigl(\theta'-1\bigr)-\mu\bigr]1_{\eta
_t(x)=1}+\bigl[-
\theta'\nu+f_0(x,\eta_t)\beta
\bigr]1_{\eta_t(x)=2}. %
\]
Noting that $\omega(\theta'-1)-\mu\le0$ and that
\[
-\theta'\nu+f_0(x,\eta_t)\beta\le-
\theta'\nu+\beta\le0 %
\]
we have shown that $\mu(\eta_t) \le0$ for all $t \ge0$. Thus,
$S(\eta_t)$ is a nonnegative supermartingale. By the martingale
convergence theorem, $S(\eta_t)$ converge to some limit as
$t\rightarrow\infty$. Note that each jump in $\eta_t$ will change
$S(\eta_t)$ by $1, \theta'$ or $\theta'-1>0$. Thus, to have
convergence of $S(\eta_t)$, with probability one there must be only
finite jumps in each path of $\eta_t$, which implies that with
probability one $\eta_t$ will end up at configuration of all 0's,
which is the absorbing state.

For the second part of the theorem, there is no nontrivial stationary
distribution when $\beta\le\nu$. When $\beta>\nu$, note that when
$\mu\nu> \omega(\beta-\nu)$, there is a $\theta'$ such that
\[
\frac{\mu+\omega}{\omega}> \theta'>\frac{\beta}{\nu}>1. %
\]
We again use the
\[
Q'(\eta_t)=\sum_{x\in\ZZ^d}
e^{-\lambda'\llVert x\rrVert }w'\bigl[\eta_t(x)\bigr] %
\]
similar to the $Q$ introduced at the beginning of Lemma \ref{reclem},
with $\lambda'>0$ and
\[
w'\bigl[\eta(x)\bigr]=\cases{ 0, &\quad if $\eta(x)=0$,
\vspace*{3pt}\cr
1, &
\quad if $\eta(x)=1$,
\vspace*{3pt}\cr
\theta', &\quad if $\eta(x)=2$.}
\]
Consider the infinitesimal mean of $Q'(\eta_t)$. Using exactly the
same argument as in Lemma \ref{infmean}, we have for any $\eta$,
\[
\mu'(\eta)\le\sum_{\eta(x)=1}\bigl[\bigl(
\theta'-1\bigr)\omega-\mu\bigr]e^{-\lambda
\llVert x\rrVert }+\sum
_{\eta(x)=2} \bigl(\beta e^{\lambda L}-\theta'\nu
\bigr)e^{-\lambda
\llVert x\rrVert }. %
\]
Thus, when $\lambda$ is small enough, $\mu'(\eta)\le0$ for all
$\eta\in\{0,1,2\}^{\ZZ^d}$ and $Q'(\eta_t)$ is a nonnegative\vspace*{1pt}
supermartingale, and thus has to converge a.s. to a limit. Then for any
$x\in\ZZ^d$, a flip at point $x$ will contribute at least
\[
e^{-\lambda\llVert x\rrVert }\min\{1,\theta-1\} %
\]
to the total value of $Q'$. So with probability one there is a
$t<\infty$ such that there is no flip at site $x$ after time $t$,
which can only correspond to the case where $\eta_s(x)\equiv0$ for
all $s\in[t,\infty)$.
\end{pf}

%s5 #&#
\section{Spatial location of the positive density box}\label{sec5}

The argument in the previous section proves the existence of a small
box $\hat{B}_x$
with positive density, but this is not useful if
we do not have control over its location. To do this, we note that the
graphical representation in Section~\ref{sec2}
shows that box process $\xi_t$ can be stochastically bounded by
Krone's model $\eta_t$ starting from the same initial configuration.
Krone's model can in turn be bounded by a branching random walk $\gamma
_t$ in which there are no deaths, 2's give birth to 2's at rate $\beta$
and births are not suppressed even if the site is occupied.

%
%le5.1 #&#
\begin{lemma} \label{BRW}\label{le5.1}
Suppose we start from $\gamma_0$ such that $\gamma_0(x)=2$ for all
$x\in\hat{B}_0$, $\gamma_0(x)=0$ otherwise.
Let $M_k(t)$ be the largest of the absolute values of the $k$th
coordinate among the occupied sites at time $t$. If $L$ is large enough
then for any $m>0$ we have
%
%e25 #&#
\begin{equation}
P\bigl( M_k(t) \ge1+(2\beta+m)L t \bigr) \le2e^{-mt}
\llvert\hat{B}_0\rrvert. \label{maxbd}
\end{equation}
From this, it follows that there is a $C_{\fontsize{8.36pt}{11pt}{\ref{BRW}}}<\infty$ so,
%
%e26 #&#
\begin{equation} \label{L2bd}
E \bigl(\bigl[M_k(t_0\log L)\bigr]^2 \bigr)
\le C_{\fontsize{8.36pt}{11pt}{\ref{BRW}}} (L\log L)^2.
\end{equation}
\end{lemma}
\begin{pf}
First, we will start from the case where $\gamma_0$ has only one
particle at~0. Rescale space by dividing by $L$. In the limit as $L\to
\infty$, we have a branching random walk $\bar\gamma_t$ with births
displaced by an amount uniform on $[-1,1]^d$. We begin by showing that
the corresponding maximum
has $E\bar M_k^2(t_0 \log L) \le C(\log L)^2$. To this, we note that
mean number of particles in $A$ at time $t$
\[
E\bigl(\bar\gamma_t(A)\bigr) = e^{\beta t} P\bigl(\bar S(t) \in
A\bigr),
\]
where $\bar S(t)$ is a random walk that makes jumps uniform on
$[-1,1]^d$ at rate $\beta$. Let $\bar S_k(t)$ be the $k$th coordinate
of $\bar S(t)$. We have
\[
E\exp\bigl(\theta\bar S_k(t)\bigr) = \exp\bigl(\beta t\bigl[\bar
\phi(\theta)- 1\bigr]\bigr)\qquad\mbox{with } \bar\phi(\theta) =
\bigl(e^\theta+ e^{-\theta}\bigr)/2. %
\]
Large deviations implies that for any $\theta>0$
\[
P\bigl( \bar S_k(t) \ge x \bigr) \le e^{-\theta x} \exp\bigl(
\beta t\bigl[\bar\phi(\theta)-1\bigr]\bigr). %
\]
By symmetry, we have that
\[
P\bigl(\bigl\llvert\bar S_k(t)\bigr\rrvert\ge x \bigr)
\le2e^{-\theta x} \exp\bigl(\beta t\bigl[\bar\phi(\theta)-1\bigr]\bigr)
\]
and hence that
%
%e27 #&#
\begin{equation}
P\bigl( \bar M_k(t) \ge x \bigr) \le2e^{-\theta x}
e^{\beta t} \exp\bigl(\beta t\bigl[\phi(\theta)-1\bigr]\bigr).
\label{LDforbar}
\end{equation}
Since the right-hand side gives the expected number of particles with
$k$th component $\ge x$.

To prove the lemma, now we return to the case $L < \infty$. Let $\phi
(\theta) = E\exp(\theta S_k(t))$
where $S(t)$ is a random walk that makes jumps uniform on $[-1,1]^d
\cap\ZZ^d/L$ at rate $\beta$.
When $\theta=1$, $\bar\phi(1)-1 = 0.543$, so if $L$ is large $\phi
(1)-1 \le1$, and by
the argument that led to (\ref{LDforbar})
\[
P\bigl( M_k(t) \ge(2\beta+m) t \bigr) \le2e^{-mt}.
\]
The\vspace*{1pt} last result is for starting for one particle at the origin. If we
start with $\llvert \hat B_0\rrvert $ particles in $\hat{B}_0/L \subset
[-1,1]^d$ in the initial configuration $\bar
\gamma_0$ then
\[
P\bigl( M_k(t) \ge1+ (2\beta+m) t \bigr) \le2\llvert
\hat{B}_0\rrvert e^{-mt}. %
\]
Taking $t=t_0\log L$, and noting that $t_0=2d/\rho>1$, gives the
desired result.
\end{pf}

%
%le5.2 #&#
\begin{lemma} \label{maxoverL}
For any $a>0$, let $M^j_k(t_0\log L)$ $1\le j \le L^a$ be the maximum
of the absolute value of $k$th coordinates in the $j$th copy of a
family of independent and identically distributed branching random walk
in Lemma \ref{BRW}. There is a $C_{\fontsize{8.36pt}{11pt}{\ref{maxoverL}}} < \infty$ so that
for large $L$
\[
P \Bigl( \max_{1\le j \le L^a} \bigl\{M^j_k(t_0
\log L)\bigr\} \ge C_{\fontsize{8.36pt}{11pt}{\ref{maxoverL}}} L \log L \Bigr) \le1/L.
\]
\end{lemma}

\begin{pf} Taking $t=t_0\log L$ in (\ref{maxbd}) and recalling
$\llvert \hat
B_0\rrvert = O(L^d)$, the right-hand side is
$\le L^{d}\exp(-mt_0$ $\log L)$ for each copy. So the probability on
the left-hand side in the lemma ${\le} L^{a+d}\exp(-mt_0\log L)$. Taking
the constant $m$ to be large enough gives the
desired result.
\end{pf}

%s6 #&#
\section{Moving particles in \texorpdfstring{$\eta_t$}{etat}}\label{sec6}

Let $H_{t,x}$ be the set of nonzero sites of $\eta_t$ in $\hat B_x$ at
time $t$.
In this section, we will use the graphical representation in Section~\ref{sec2} and
an argument from Durrett and Lanchier \cite{DL0} to show that

%
%le6.1 #&#
\begin{lemma}\label{MoS}
There are constants $\delta_{\fontsize{8.36pt}{11pt}{\ref{MoS}}}>0 $
and an $L_{\fontsize{8.36pt}{11pt}{\ref{MoS}}} <\infty$ such that for all $L>L_{\fontsize{8.36pt}{11pt}{\ref{MoS}}}$
and any initial configuration $\eta_0$ with $\llvert H_{0,0}\rrvert
\ge L^{d/2}$
\[
P\bigl(\llvert H_{1,v}\rrvert\ge\delta_{\fontsize{8.36pt}{11pt}{\ref{MoS}}} \llvert
H_{0,0}\rrvert\bigr) > 1-e^{-L^{d/4}} %
\]
for any $v \in\{0, \pm e_1,\ldots,\pm e_d\}$.
\end{lemma}

\begin{pf}
We begin with the case $v=0$ which is easy. Define $G^0_0$ to be the
set of points $x\in\hat{B}_0$,
with (a) $\eta_0(x)\ge1$, and (b) no death marks $\times$'s occur in
$\{x\}\times[0,1]$.
We have $\xi_t(x)\ge1$ on $S_0 = H_{0,0} \cap G_0^0$, and $\llvert
S_0\rrvert \sim
\operatorname{Binomial}(\llvert H_{0,0}\rrvert,e^{-\mu})$,
so the desired result follows from large deviations for the Binomial.

For $v\neq0$, define $G_0$ to be the set of points in $G^0_0$ for
which (c) there exists a $(\bullet)$,
which produces growth from type 1 to type 2, in $\{x\}\times[0,1/2]$.
We define $G_v$ to be the set of points $y$ in $\hat{B}_v$ so that
there are no $\times$'s in $\{y\}\times[0,1]$.
For any $x\in\hat{B}_0$ and $y\in\hat{B}_v$ we say that $x$ and $y$
are connected (and write $x\to y$) if there is an
arrow from $x$ to $y$ in $(1/2,1)$. By definition of our process $\eta
_1(y)\ge1$ for all $y$ in
\[
S=\{y\dvtx y\in G_v\mbox{, there exists an $x\in G_0$ so
that $x\to y$} \}. %
\]
It is easy to see that
%
%e28 #&#
\begin{equation}
\llvert G_0\rrvert\sim \operatorname{Binomial}\bigl[\llvert H_{0,0}
\rrvert,e^{-\mu}\bigl(1-e^{-\omega/2}\bigr)\bigr]. \label{DoGS}
\end{equation}
Conditional on $\llvert G_0\rrvert $:
%
%
%e29 #&#
\begin{equation}
\label{CdoS} \llvert S\rrvert\sim \operatorname{Binomial}\bigl(\llvert\hat B_v
\rrvert,e^{-\mu}\bigl[1-e^{-\beta\llvert G_0\rrvert /2\llvert
B_0\rrvert }\bigr]\bigr),
\end{equation}
by Poisson thinning since the events of being the recipient of a birth
from $\hat B_0$
are independent for different sites in $G_v$.

Since the binomial distribution decays exponentially fast away from the
mean, there is some constant $c>0$ such that
%
%e30 #&#
\begin{equation}
P\bigl(\llvert G_0\rrvert>\llvert H_{0,0}\rrvert
e^{-\mu}\bigl(1-e^{-\omega/2}\bigr)/2\bigr)\ge1-e^{-cL^{d/2}}.
\label{A1}
\end{equation}
To simplify the next computation, we note that $1-e^{-\beta r} \sim
\beta r$ as $r\to0$ so if the $\varepsilon_0$
in the definition of the small box is small enough
\[
1-e^{-\beta\llvert G_v\rrvert /2\llvert B_0\rrvert } \ge\beta
\llvert G_v\rrvert/4\llvert
B_0\rrvert. %
\]
Let $p=e^{-\mu}\beta\llvert G_0\rrvert /4\llvert B_0\rrvert $.
A standard large deviations result, see, for example, Lemma 2.8.5 in
\cite
{RGD} shows that
if $X=\operatorname{Binomial}(N,p)$ then
\[
P(X \le Np/2) \le\exp(-Np/8) %
\]
from which the desired result follows.
\end{pf}

Let $\llVert \cdot\rrVert _1$ be the $L^1$-norm on $\ZZ^d$.
Our next step is to use Lemma \ref{MoS} $O(\log L)$ times to prove:

%
%le6.2 #&#
\begin{lemma}\label{MNS}
For any $\alpha\in(d/2,d)$, let $C_{\fontsize{8.36pt}{11pt}{\ref{MNS}}}$ be a constant such
that $C_{\fontsize{8.36pt}{11pt}{\ref{MNS}}} \log\delta_{\fontsize{8.36pt}{11pt}{\ref{MoS}}}>\alpha-d$.
There is a finite $L_{\fontsize{8.36pt}{11pt}{\ref{MNS}}} > L_{\fontsize{8.36pt}{11pt}{\ref{MoS}}}$ such that for all
$L>L_{\fontsize{8.36pt}{11pt}{\ref{MNS}}}$,
any initial configuration $\eta_0$ with $\llvert H_{0,0}\rrvert \ge
a_0\llvert \hat{B}_0\rrvert $,
and any $x\in\ZZ^d$ such that $\llVert x\rrVert _1\le C_{\fontsize{8.36pt}{11pt}{\ref{MNS}}}
\log L$, we have
\[
P \bigl(\llvert H_{x,[C_{\fontsize{8.36pt}{11pt}{\ref{MNS}}} \log L]}\rrvert\ge L^{\alpha
} \bigr)
\ge1-e^{-L^{d/4}/2} .
\]
\end{lemma}
\begin{pf} Let $n=[C_{\fontsize{8.36pt}{11pt}{\ref{MNS}}} \log L]$. We can find a sequence
$x_0=0,x_1,\ldots, x_n=x$ such that for all $i=0,\ldots,n-1$,
$x_{i-1}-x_i\in\{0, \pm e_1,\ldots,\pm e_d\}$. For any $i=1,\ldots,n$
define the event
\[
A_i= \bigl\{\llvert H_{x_i,i}\rrvert\ge
\delta_{\fontsize{8.36pt}{11pt}{\ref{MoS}}}^i \llvert H_{0,0}\rrvert\bigr\}.
\]
By the definition of $C_{\fontsize{8.36pt}{11pt}{\ref{MNS}}}$, $\llvert H_{x,[C_{\fontsize{8.36pt}{11pt}{\ref{MNS}}}
\log L]}\rrvert \ge L^{\alpha}$ on $A_n$.
To estimate $P(A_n)$ note that by Lemma \ref{MoS}
\begin{eqnarray*}
P(A_n) &\ge&1-\sum_{i=1}^n P
\bigl(A_i^c\bigr)\ge1-\sum_{i=1}^n
P\bigl(A_i^c\mid A_{i-1}\bigr)
\\
& \ge&1-C_{\fontsize{8.36pt}{11pt}{\ref{MNS}}}( \log L ) e^{-L^{d/4}}\ge1-e^{-L^{d/4}/2}
\end{eqnarray*}
when $L$ is large.
\end{pf}

%s7 #&#
\section{Block construction and the proof of Theorem~\texorpdfstring{\protect\ref{teo1}}{1}}\label{sec7}

At this point, we have all the tools to construct the block event and
complete the proof of Theorem~\ref{teo1}. Let $0 < a < \alpha/2-d/4$,
$K=L^{1+2a/3}$, $\Gamma_m = 2mKe_1 + [-K,K]^d$, and $\Gamma'_m =
2mKe_1 + [-K/2,K/2]^d$. If $m+n$ is even, we say that $(m,n)$ is wet if
there is a positive\vadjust{\goodbreak} density small box, that is, a box with size $\ell$
and densities of nonzero sites $\ge a_0$, in~$\Gamma'_m$ at some time in $[nL^a,nL^a+C_6\log L]$, where $C_6 =
C_{\fontsize{8.36pt}{11pt}{\ref{MNS}}}+t_0$. Our goal is to show

%le7.1 #&#
\begin{lemma}\label{block}
If $(m,n)$ is wet then with high probability so is $(m+1,n+1)$, and the
events which produce this are measurable with respect to the graphical
representation in
$(\Gamma_m \cup\Gamma_{m+1}) \times[nL^a,(n+1)L^a + C_6\log L]$.
\end{lemma}

Once this is done, Theorem~\ref{teo1} follows. See \cite{DStF} for more details.

\begin{pf*}{Proof of Lemma \ref{block}}
To prove Lemma \ref{block}, we will alternate two steps, starting from
the location of the initial positive density box $\hat B_{y_0}$ at time $T_0$.
Let $A_0 = \{T_0 < \infty\}$ which is the whole space. Assume given a
deterministic sequence $\delta_i$
with $\llVert \delta_i\rrVert \le C_{\fontsize{8.36pt}{11pt}{\ref{MNS}}}\log L$. If we never
meet a
failure, the construction will terminate at the first time that $T_i >
(n+1)L^a$.
The actual number steps will be random but the number is ${\le} N =
\lceil L^a/[C_{\fontsize{8.36pt}{11pt}{\ref{MNS}}} \log L] \rceil$.
We will estimate the probability of success supposing that $N$ steps
are required.
This lower bounds the probability of success when we stop at the first
time $T_i \ge(n+1)L^a$. Suppose $i \ge1$.

\begin{longlist}
\item[\textit{Deterministic moving}.] If at the stopping time $T_{i-1}<\infty$,
we have a positive density small box $\hat{B}_{y_{i-1}}$, then we
use results in Section~\ref{sec6} to produce a small box $\hat
{B}_{y_{i-1}+\delta_i}$ with at least $L^{\alpha}$ nonzero sites at
time $S_i=T_{i-1}+[C_{\fontsize{8.36pt}{11pt}{\ref{MNS}}}\log L]$.
If we fail, we let $S_i=\infty$ and the construction terminates. Let
$A_i^+ = \{S_i < \infty\}$.
\end{longlist}

\begin{longlist}
\item[\textit{Random recovery}.] If at the stopping time $S_i < \infty$, we have
a small box $\hat{B}_{y_{i-1}+\delta_i}$ with at least $L^\alpha$
nonzero sites then we set all of the sites outside the box to 0, and we
use the recovery lemma to produce
a positive density small box $\hat B_{y_i}$ at time $S_i \le T_i \le
S_i + t_0 \log L$. Again if we fail,
we let $T_i=\infty$ and the construction terminates. Let $A_i = \{T_i
< \infty\}$.
Let $\Delta_i(\omega) = y_i - (y_{i-1}+\delta_i)$ on $A_i$, and $=0$
on~$A_i^c$.
\end{longlist}

If we define the partial sums $\bar y_i = y_0 + \sum_{j=1}^i\delta_i$
and $\Sigma_i=\sum_{j=1}^i \Delta_j$,
then we have $y_i= \bar y_i + \Sigma_i$. We think of $\bar y_i$ as the mean
of the location of the positive density box and $\Sigma_i$ as the
random fluctuations in its location.
We make no attempt to adjust the deterministic movements $\delta_i$ to
compensate for the fluctuations. Let $y_{\mathrm{end}}=(y_{\mathrm{end}}^1,0,\ldots,
0)\in\ZZ^d$ be such that $2K(m+1)e_1\in\hat{B}_{y_{\mathrm{end}}}$. We
define the $\delta_i$
to reduce the coordinates $y_0^k$, $k=2, \ldots, d$ to 0 and then
increase $y_0^1$ to $y^1_{\mathrm{end}}$,
in all cases using steps of size $\le C_{\fontsize{8.36pt}{11pt}{\ref{MNS}}} \log L$. Note that
\[
\llVert y_0-y_{\mathrm{end}}\rrVert_1=O(K)/\ell=O
\bigl(L^{2a/3}\bigr)=o(N)
\]
so we can finish the movements well before $N$ steps. And once this is
done we set the remaining $\delta_i$ to 0. Moreover, note that each
successful step in our iteration takes a time at most $C_6\log L$.
Thus, we will get to $y_{\mathrm{end}}$ by $t=nL^a+O(L^{2a/3})C_6\log
L<(n+1)L^a$. At the first time $T_i\ge(n+1)L^a$, we already have
%
%e31 #&#
\begin{equation}
\label{RbT} \delta_i=0,\qquad \bar{y}_i=y_{\mathrm{end}}.
\end{equation}

At this point, we are ready to state the main lemma of this section
that controls the spatial movement in our iteration.

%
%le7.2 #&#
\begin{lemma}\label{MaP}
For any initial configuration $\eta(T_0)$ so that there is a small box
$\hat{B}_{y_0}\subset\Gamma'_m$,
and any sequence $\delta_i$, $i\le N$ with $\llVert \delta_i\rrVert
\le C_{\fontsize{8.36pt}{11pt}{\ref{MNS}}} \log L$ and any $\varepsilon>0$, there is a good event $G_N$
with $G_N \to1$ as $L \to\infty$ so that \textup{(a)} $G_N \subset A_N$, \textup{(b)}
on $G_N$, $\llVert y_i-\bar y_i\rrVert <\varepsilon L^{2a/3}$ for
$1\le i \le N$,
\textup{(c)} $G_N$ depends only on the gadgets of graphical representation in
$\Gamma_m \cup\Gamma_{m+1}$.
\end{lemma}

\begin{pf} The first step is to show that $P(A_N) \to1$ as $L\to
\infty$.
For the $i$th deterministic moving step, using the strong Markov
property and Lemma \ref{MNS}, we have
\[
P\bigl(A_{i}^+\mid A_{i-1}\bigr)>1-e^{-L^{d/4}/2}.
\]
Then for the random recovery phase, according to Lemma \ref{reclem},
we have the conditional probability of success:
\[
P\bigl(A_i\mid A_{i}^+\bigr)=P(\tau_i<t_0
\log L)>1-L^{d/2-\alpha}. %
\]
Combining the two observations, we have
%
%
%e32 #&#
\begin{equation}
\label{PoS} P(A_i\mid A_{i-1})>\bigl(1-e^{-L^{d/4}/2}
\bigr) \bigl(1-L^{d/2-\alpha
}\bigr)>1-e^{-L^{d/4}/2}-L^{d/2-\alpha}
\end{equation}
which implies
\begin{eqnarray*}
P(A_N) &\ge& 71-\sum_{i=1}^N
P\bigl(A_i^c\bigr) \ge1-\sum
_{i=1}^NP\bigl(A_i^c\mid
A_{i-1}\bigr)
\\
& \ge&1-L^a\bigl(e^{-L^{d/4}/2}+L^{d/2-\alpha}\bigr)
\ge1-2L^{d/4-\alpha/2} \to1.
\end{eqnarray*}

The next step is to control the fluctuations in the movement of our box.

%
%le7.3 #&#
\begin{lemma}\label{MaBV}
Let ${\mathcal F}(T_i)$ be the filtration generated by events in the
graphical representation up to stopping time $T_i$.
For any $1 \le k \le d$, $\{\Sigma_i^k\}_{i=1}^{N}$ is a martingale
with respect to $\mathcal{F}(T_{i})$.
$E(\Sigma_i^k)=0$, and $\var(\Sigma_{N}^k)\le C_{\fontsize{8.36pt}{11pt}{\ref{MaBV}}}
L^{a}\log L$ so for any $\varepsilon$ we have
\[
P \Bigl( \max_{i \le N} \llVert\Sigma_i \rrVert>
\varepsilon L^{2a/3} \Bigr) \to0 .
\]
\end{lemma}

\begin{pf} Consider the conditional expectation of $\Delta_i^k$
under $\mathcal{F}(S_i)$.
According to the discussions about the truncated process right before
Section~\ref{sec2.1}, we have $E(\Delta_i^k\mid{\mathcal F}(S_i))=0$. Noting that
$2\ell\llvert \Delta_i^k\rrvert $ can be bounded by the largest $k$th
coordinate
among the occupies sites of the corresponding branching random walk at
time $t_0\log L$, Lemma \ref{BRW} implies that
$E((\Delta_i^k)^2\mid\break {\mathcal F}(S_i))\le C(\log L)^2$. By
orthogonality of
martingale increments
$\var(\Sigma_{N}^k)\le N C(\log L)^2$. Since $N \le L^a/[C_{5.2} \log
L]+1$, we have the desired bound on the variances and
the desired result follows from $L^2$ maximal inequality for martingales.\vadjust{\goodbreak}
\end{pf}

To check (c), now note that under $A_{i}$ the success of $A_{i+1}^+$
depends only on gadgets in
\[
\bigl(2\ell y_i+[-\ell C_{\fontsize{8.36pt}{11pt}{\ref{MNS}}}\log L,\ell
C_{\fontsize{8.36pt}{11pt}{\ref{MNS}}}\log L]^d \bigr) \times[T_i,S_{i+1}]
\]
and that when the $i$th copy of the truncated process never wanders outside
\[
D_i=2\ell(y_{i-1}+ \delta_i) +
\bigl[-C_{\fontsize{8.36pt}{11pt}{\ref{maxoverL}}}(L \log L), C_{\fontsize{8.36pt}{11pt}{\ref{maxoverL}}}(L \log L)
\bigr]^d %
\]
the success of $A_{i+1}$ under $A^+_{i+1}$ depends only on gadgets in
$D_i$. According to Lemma \ref{MaP} and the fact that $N<L^a$, with
probability $\ge1-L^{-1}=1-o(1)$, our construction only depends on
gadgets in the box:
\begin{eqnarray*}
&& \bigcup_{i=0}^{N-1}  \bigl[\bigl(2\ell
y_i+[-\ell C_{\fontsize{8.36pt}{11pt}{\ref{MNS}}}\log L,\ell C_{\fontsize{8.36pt}{11pt}{\ref{MNS}}}\log
L]^d \bigr) \times[T_i,S_{i+1}]
\\
&&\hspace*{20pt}{} \cup
\bigl(2\ell(y_{i}+ \delta_{i+1}) + \bigl[-C_{\fontsize{8.36pt}{11pt}{\ref{maxoverL}}}(L
\log L), C_{\fontsize{8.36pt}{11pt}{\ref{maxoverL}}}(L \log L)\bigr]^d\bigr)
\times[S_{i+1},T_{i+1}] \bigr].
\end{eqnarray*}
The locations of the $y_i$ are controlled by Lemma \ref{MaBV} so that
it is easy to see that the box defined above is a subset of $\Gamma
_n\cup\Gamma_{n+1}$, and proof of Lemma \ref{MaP} is complete.
\end{pf}

Back to the proof of Lemma \ref{block}, on $G_N$, it follows from
(\ref{RbT}) and Lemma \ref{MaP}, when we stop at the first time $T_i
\ge(n+1)L^a$:
\[
\bigl\llVert y_i-2K(m+1)e_1\bigr\rrVert\le\ell
\bigl(1+2\llVert y_i-\bar{y}_i\rrVert\bigr)\le4\varepsilon
K %
\]
which implies
\[
\hat{B}_{y_i}=2\ell y_i+(-\ell,\ell]^d\subset
\Gamma_{n+1}'. %
\]
Noting that the success of $G_N$ only depends on gadgets in $\Gamma
_n\cup\Gamma_{n+1}$, we have proved that $G_N$ is measurable with
respect to the space--time box in the statement of Lemma \ref{block},
which completes the proof of Lemma \ref{block} and Theorem \ref{mainth}.
\end{pf*}

%s8 #&#
\section{Proof of Theorem \texorpdfstring{\protect\ref{SurHD}}{3}}\label{sec8}

\begin{figure}%[b]

\includegraphics{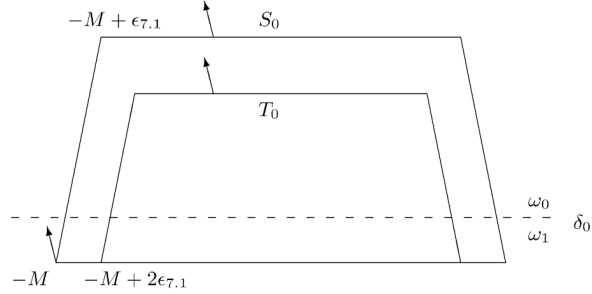}

%%
%\begin{center}
%%
%\begin{picture}(300,130)
%\put(50,10){\line(1,0){200}}
%\put(70,110){\line(1,0){160}}
%\put(50,10){\line(1,5){20}}
%\put(230,110){\line(1,-5){20}}
%\put(70,10){\line(1,5){15}}
%\put(85,85){\line(1,0){130}}
%\put(215,85){\line(1,-5){15}}
%\put(30,30){\dashline{4}(0,0)(240,0)}
%\put(280,25){$\delta_0$}
%\put(260,35){$\omega_0$}
%\put(260,20){$\omega_1$}
%\put(140,115){$S_{0}$}
%\put(140,75){$T_{0}$}
%\put(30,0){$-M$}
%\put(55,115){$-M+\varepsilon_{7.1}$}
%\put(62,0){$-M+2\varepsilon_{7.1}$}
%\put(50,10){\vector(-1,4){4}}
%\put(120,110){\vector(-1,4){4}}
%\put(120,85){\vector(-1,4){4}}
%\end{picture}
%%
%\end{center}
%%
\caption{Test functions for $d=1$.}\label{fig1}
\end{figure}

Our \textit{first step} is to construct the test functions for $S$ and $T$
and show that, under (\ref{IDE}), they have positive derivatives for
all sites in the region of interest. According to (\ref{SFsur2}), we
can choose $\Sigma_0\in(0,1)$ such that
%
%e33 #&#
\begin{equation}
\label{STsum} \frac{\nu}{\omega_0}<\frac{\beta(1-\Sigma_0)}{2^{d}(\mu
+\omega_0)}.
\end{equation}
Let
%
%e34 #&#
\begin{equation}
\label{STratio} \gamma_0\in\biggl(\frac{\nu}{\omega_0},
\frac{\beta(1-\Sigma
_0)}{2^{d}(\mu+\omega_0)} \biggr)
\end{equation}
and let
%
%e35 #&#
\begin{equation}
\label{DoB} T_0=\frac{\Sigma_0}{1+\gamma_0}, \qquad S_0=
\frac{\gamma_0\Sigma
_0}{1+\gamma_0}.
\end{equation}
Note that
%
%e36 #&#
\begin{equation}
\label{SnR} S_0+T_0=\Sigma_0,\qquad
S_0/T_0=\gamma_0.
\end{equation}
Recall that for any $x\in\ZZ^d$ and $r\ge0$, $B(x,r)$ is defined in
Section~\ref{sec2} to be the $L^\infty$ neighborhood of $x$ with range $r$.
With $S_0, T_0$ defined as above and $\varepsilon_{\fontsize{8.36pt}{11pt}{\ref{DIDE}}}$ to be
specified later, define the test functions $S_{\mathrm{test}}(x,0)$ and
$T_{\mathrm{test}}(x,0)$ as follows (Figure~\ref{fig1}
%
%f1 #&#
shows those test
functions when
$d=1$): let $S_{\mathrm{test}}(x,0)=S_0$ on $B(0,M-\varepsilon_{\fontsize{8.36pt}{11pt}{\ref{DIDE}}})$,
$S_{\mathrm{test}}(x,0)=0$ on $ B(0,M)^c$, and
%
%
%e37 #&#
\begin{equation}
\label{TestS} S_{\mathrm{test}}(x,0)=\frac
{d[x,B(0,M)^c]S_0}{d[x,B(0,M)^c]+d[x,B(0,M-\varepsilon
_{\fontsize{8.36pt}{11pt}{\ref{DIDE}}})]}
\end{equation}
for $x\in B(0,M-\varepsilon_{\fontsize{8.36pt}{11pt}{\ref{DIDE}}})^c\cap B(0,M)$. Similarly, let
$T_{\mathrm{test}}(x,0)=T_0$ on $B(0,M-3\varepsilon_{\fontsize{8.36pt}{11pt}{\ref{DIDE}}})$,
$T_{\mathrm{test}}(x,0)=0$ on $B(0,M-2\varepsilon_{\fontsize{8.36pt}{11pt}{\ref{DIDE}}})^c$, and
%
%
%e38 #&#
\begin{equation}
\label{TestT} T_{\mathrm{test}}(x,0)=\frac{d[x,B(0,M-2\varepsilon_{\fontsize{8.36pt}{11pt}{\ref{DIDE}}})^c]T_0}{d[x,B(0,M-2\varepsilon_{\fontsize{8.36pt}{11pt}{\ref{DIDE}}})
^c]+d[x,B(0,M-3\varepsilon_{\fontsize{8.36pt}{11pt}{\ref{DIDE}}})]}
\end{equation}
for $ x\in B(0,M-3\varepsilon_{\fontsize{8.36pt}{11pt}{\ref{DIDE}}})^c\cap B(0,M-2\varepsilon
_{\fontsize{8.36pt}{11pt}{\ref{DIDE}}})$. In the definitions above, $M=\max\{4,4\kappa\}$, $d(x,A)$
be the distance between $x\in R^d$ and $A\subset R^d$ under $L^\infty
$-norm, $\varepsilon_{\fontsize{8.36pt}{11pt}{\ref{DIDE}}}=\varepsilon_{\fontsize{8.36pt}{11pt}{\ref{DIDE}}}(\beta,\omega
_0,\mu,\nu,\kappa)$ is some positive constant that will be specified
later in the proof of Lemma \ref{DIDE}. The following lemma shows that
the test functions have positive derivatives under IDE (\ref{IDE}).

%re2 #&#
\begin{rem}
Throughout the discussion in this section, all the
$\varepsilon$'s, $\delta$'s, $t$'s and $c$'s introduced are constants
independent to the choice of $L$.
\end{rem}

%
%le8.1 #&#
\begin{lemma}\label{DIDE}
Under the conditions in Theorem \ref{SurHD}, there are $\varepsilon
_{\fontsize{8.36pt}{11pt}{\ref{DIDE}}}$ and $\varepsilon_1>0$ so that under IDE (\ref{IDE}):
%
%e39 #&#
\begin{eqnarray}
\label{Pdift} \frac{dS_{\mathrm{test}}(x,0)}{dt}&\ge&4\varepsilon_1\qquad\mbox{for all } x
\in B(0,M)=\overline{\{S_{\mathrm{test}}>0\}},
\nonumber\\[-8pt]\\[-8pt]\nonumber
\frac{dT_{\mathrm{test}}(x,0)}{dt}&\ge&4\varepsilon_1\qquad\mbox{for all } x\in B(0,M-2
\varepsilon_{\fontsize{8.36pt}{11pt}{\ref{DIDE}}})=\overline{\{T_{\mathrm{test}}>0\}}.
\end{eqnarray}
\end{lemma}

\begin{pf}
With $T_{\mathrm{test}}$ and $S_{\mathrm{test}}$ defined as above, for any $x\in B(0,M)$,
the local grass density can be upper bounded as follows (see
Figure~\ref{STr} for the case when $d=2$):
\[
D_\kappa^G(x,0)\le1- \biggl(\frac{\kappa-\varepsilon_{\fontsize{8.36pt}{11pt}{\ref{DIDE}}}}{2\kappa}
\biggr)^d S_0\le1-2^{-d} \biggl(1-d
\frac{\varepsilon
_{\fontsize{8.36pt}{11pt}{\ref{DIDE}}}}{\kappa} \biggr)S_0. %
\]
%
%
%f2 #&#
\begin{figure}%[b]

\includegraphics{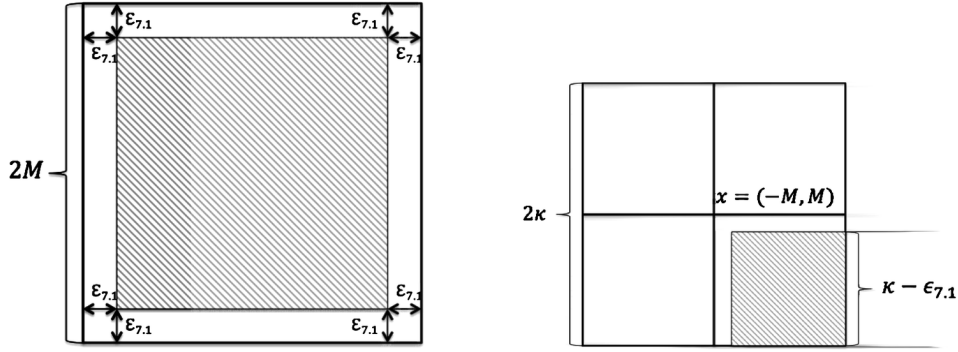}

\caption{$S_{\mathrm{test}}$ with $d=2$. Left: $S_{\mathrm{test}}$ with $d=2$. The big
box of size $2M$ is region $\overline{\{S_{\mathrm{test}}>0\}}$. The shadow
area is the region where $S_{\mathrm{test}}=S_0$. Right: Worst case for
$S_{\mathrm{test}}$ where $x=(-M,M)$. The big box is the local grass
environment. The shadow area is the region where $S_{\mathrm{test}}=S_0$.}
\label{STr}
\end{figure}

Noting that $\delta_0<2^{-d}S_0$, let
%
%e40 #&#
\begin{equation}
\label{Tshift1} \varepsilon_{\fontsize{8.36pt}{11pt}{(\ref{Tshift1})}}=\frac{\kappa(1-2^d\delta
_0/S_0)}{4d}.
\end{equation}
It is easy to check that when $\varepsilon_{\fontsize{8.36pt}{11pt}{\ref{DIDE}}}\le\varepsilon
_{\fontsize{8.36pt}{11pt}{(\ref{Tshift1})}}$
\[
D_\kappa^G(x,0)< 1-\delta_0, %
\]
which implies that
%
%e41 #&#
\begin{equation}
\label{PGR} \omega\bigl[D_\kappa^G(x,0)\bigr]\equiv
\omega_0
\end{equation}
for all $x\in B(0,M)$. That is, all sites in the region of test function
live in a environment with a higher growth rate $\omega_0$. It is easy
to see that the derivative of $T_{\mathrm{test}}$ can be lower bounded by its
derivative on the top, that is, for any $x\in B(0,M-2\varepsilon_{\fontsize{8.36pt}{11pt}{\ref{DIDE}}})$
%
%
%e42 #&#
\begin{equation}
\label{DTT} \frac{dT_{\mathrm{test}}(x,0)}{dt}\ge\frac{dT_{\mathrm{test}}(0,0)}{dt}\ge
\omega
_0S_{0}-\nu T_{0}
\end{equation}
and this holds for all $\kappa>0$. Note that $\gamma_0=S_0/T_0$
according to (\ref{DoB}). Combining this observation with the
definition of $\gamma_0$ in (\ref{STratio})
%
%e43 #&#
\begin{equation}
\label{PDT} \omega_0S_{0}-\nu T_{0}=
\omega_0T_0 \biggl(\frac{S_0}{T_0}-\frac{\nu
}{\omega_0}
\biggr)=\omega_0T_0 \biggl(\gamma_0-
\frac{\nu}{\omega
_0} \biggr)>0.
\end{equation}
Thus, we have the derivative of $T_{\mathrm{test}}$ is always positive for all
$x\in B(0,M-2\varepsilon_{\fontsize{8.36pt}{11pt}{\ref{DIDE}}})$.

Similarly, we can control the lower bound of derivative for test
function $S_{\mathrm{test}}$ as follows: for any $x\in B(x,M)$
%
%
%e44 #&#
\begin{eqnarray}
\label{DST} \frac{dS_{\mathrm{test}}(x,0)}{dt}&\ge&\frac{\beta
T_{0}(1-3\varepsilon_{\fontsize{8.36pt}{11pt}{\ref{DIDE}}})^d}{2^d}(1-S_{0}-T_0)-(
\mu+\omega_0)S_{0}
\nonumber\\[-8pt]\\[-8pt]\nonumber
&\ge&\frac{\beta T_{0}(1-3d\varepsilon_{\fontsize{8.36pt}{11pt}{\ref{DIDE}}})}{2^d}(1-S_{0}-T_0)-(\mu+
\omega_0)S_{0}.
\end{eqnarray}
For the right-hand side of (\ref{DST}), according to (\ref{SnR})
\begin{eqnarray*}
\frac{dS_{\mathrm{test}}(x,0)}{dt}&\ge& \frac{\beta T_{0}(1-3d\varepsilon_{\fontsize{8.36pt}{11pt}{\ref{DIDE}}})}{2^d}(1-S_{0}-T_0)-(
\mu+\omega_0)S_{0}
\\
&=&T_0(\mu+\omega_0) \biggl[(1-3d\varepsilon_{\fontsize{8.36pt}{11pt}{\ref{DIDE}}})
\frac{\beta
(1-S_0-T_0)}{2^{d}(\mu+\omega_0)}-\frac{S_0}{T_0} \biggr]
\\
&=&T_0(\mu+\omega_0) \biggl[(1-3d\varepsilon_{\fontsize{8.36pt}{11pt}{\ref{DIDE}}})
\frac{\beta
(1-\Sigma_0)}{2^{d}(\mu+\omega_0)}-\gamma_0 \biggr].
\end{eqnarray*}
Again recalling the definition in (\ref{STratio}) that
\[
\frac{\beta(1-\Sigma_0)}{2^{d}(\mu+\omega_0)}>\gamma_0, %
\]
we let
%
%e45 #&#
\begin{equation}
\label{Tshift2} \varepsilon_{\fontsize{8.36pt}{11pt}{(\ref{Tshift2})}}= \biggl[ 1- \frac
{2^{d}(\mu+\omega
_0)\gamma_0} {\beta(1-\Sigma_0)}
\biggr]\Big/(6d).
\end{equation}
So for any $\varepsilon_{\fontsize{8.36pt}{11pt}{\ref{DIDE}}}\le\varepsilon_{\fontsize{8.36pt}{11pt}{(\ref{Tshift2})}}$
and $x\in B(0,M)$
%
%e46 #&#
\begin{equation}
\label{PDS} \frac{dS_{\mathrm{test}}(x,0)}{dt}\ge T_0(\mu+\omega_0)
\biggl[(1-3d\varepsilon_{\fontsize{8.36pt}{11pt}{(\ref{Tshift2})}})\frac{\beta(1-\Sigma
_0)}{2^{d}(\mu+\omega
_0)}-\gamma_0
\biggr]>0.
\end{equation}
Thus, let
%
%e47 #&#
\begin{equation}
\label{EIDE} \varepsilon_{\fontsize{8.36pt}{11pt}{\ref{DIDE}}}=\min\{\varepsilon_{\fontsize{8.36pt}{11pt}{(\ref
{Tshift1})}},\varepsilon
_{\fontsize{8.36pt}{11pt}{(\ref{Tshift2})}}\}>0
\end{equation}
and
\begin{eqnarray*}
\varepsilon_1&=&\frac{1}{4}\min \biggl\{\omega_0T_0
\biggl(\gamma_0-\frac{\nu}{\omega_0} \biggr), T_0(\mu+
\omega_0) \biggl[(1-3d\varepsilon_{\fontsize{8.36pt}{11pt}{(\ref{Tshift2})}})\frac{\beta(1-\Sigma_0)}{2^{d}(\mu+\omega_0)}-
\gamma_0 \biggr] \biggr\}
\\
&>&0.
\end{eqnarray*}
Combining (\ref{PDT}) and (\ref{PDS}) and the proof is complete.
\end{pf}

With the test functions constructed, our \textit{second step} is similar
to the proof of Theorem \ref{mainth}. We introduce the truncated
version of the Staver--Levin model, and as before, denote the process
by $\bar\xi_t$. For $\ell=\varepsilon L$, where the exact value of
$\varepsilon$ is specified later in Lemma \ref{ExRe}, $\bar\xi_t$ has
birth rate $\beta\llvert \mathcal{N}(x,L)\rrvert /(2L+1)^d$, where
$\mathcal
{N}(x,L)$ in the truncated neighborhood defined in Section~\ref{sec2}. Type 1's
and 2's in $\bar\xi_t$ die at the same rates as in the original $\chi
_t$, while a growth of a sapling into a tree occurs at rate:
\[
\bar\omega\bigl[\bar G(x,\xi)\bigr]=\cases{ \omega_0, &\quad$\bar
G(x,\xi)\in[0,1-\delta)$,
\vspace*{3pt}\cr
\omega_1, &\quad$\bar G(x,\xi)\in[1-
\delta,1]$,}
\]
where
\[
\bar G(x,\xi)=\frac{\mbox{\# of 0's in }\mathcal
{N}(x,K)}{\llvert \mathcal{N}(x,K)\rrvert }, %
\]
$K=\kappa L$, and $\delta=\delta_0+4d\varepsilon$. First of all, with
the same argument as in Section~\ref{sec2}, we immediately have that the number
of different types in each small box forms a Markov process $\bar\zeta
_t$. According to (\ref{NxLB}):
\begin{eqnarray*}
\bar G(x,\xi)&\le&1-\frac{\mbox{\# of $(1+2)$'s in }\mathcal
{N}(x,K)}{\llvert B(x,K)\rrvert }
\\
&\le& G(x,\xi)+1-(1-4\varepsilon)^d\le G(x,\xi)+4d\varepsilon,
\end{eqnarray*}
combining this with the definition of $\delta$, we have for any $\xi
'\ge\xi$, $\omega[G(x,\xi')]\ge\bar\omega[\bar G(x,\xi)]$,
which implies that the truncated $\bar\xi_t$ once again is dominated
by the original $\chi_t$. Thus, in order to prove Theorem \ref
{SurHD}, it suffices to show that $\bar\xi_t$ survives.

The \textit{third step} is to construct a initial configuration of $\bar
\xi_0$ according to the test functions defined in (\ref{TestS}) and
(\ref{TestT}). For any $x\in\ZZ^d$: if $2\ell x\notin B(0,ML)$ there
is no saplings\vspace*{1pt} or trees in $\hat{B}_x$ under $\bar\xi_0$. If $2\ell
x\in B(0,ML)$, $n_1(x,\bar\xi_0)=\llvert \hat B_0\rrvert
S_{\mathrm{test}}(2\ell x/L)$,
$n_2(x,\bar\xi_0)=\llvert \hat B_0\rrvert T_{\mathrm{test}}(2\ell x/L)$. As
noted earlier in
the box process $\bar\xi_t$, the locations of the 1's and 2's inside
each small box makes no difference.

We then look at $f_i(x,\bar\xi_t)$, the densities of type $i$ in each
small box. For any $x$, the infinitesimal means of $f_1$ and $f_2$ can
be written as follows:
%
%
%e48 #&#
\begin{eqnarray}
\label{InfM12} \mu_1(x,\bar\xi)&=&(2L+1)^{-d}\llvert\hat
B_0\rrvert\biggl(\sum_{y\dvtx \hat
B_y\subset\mathcal{N}(x,L)}
f_2(y,\bar\xi) \biggr)f_0(x,\bar\xi)\beta\nonumber
\\
&&{} -\bigl[\bar\omega\bigl(\bar G(x,\bar\xi)\bigr)+\mu\bigr]f_1(x,
\bar\xi),
\\
\mu_2(x,\bar\xi)&=&\bar\omega\bigl(\bar G(x,\bar\xi)
\bigr)f_1(x,\bar\xi)-\nu f_2(x,\bar\xi).\nonumber
\end{eqnarray}
We prove the following.

%
%le8.2 #&#
\begin{lemma}\label{PIM}
There is a $\varepsilon_{\fontsize{8.36pt}{11pt}{\ref{PIM}}}>0$ such that for any $\ell
=\varepsilon
L\le\varepsilon_{\fontsize{8.36pt}{11pt}{\ref{PIM}}}L$ and the $\bar\xi_0$ defined above, we have:
\begin{itemize}
\item$\mu_1(x,\bar\xi_0)\ge2\varepsilon_{1}$ for all $2\ell x\in
B(0,ML+4\ell)$.
\item$\mu_2(x,\bar\xi_0)\ge2\varepsilon_{1}$ for all $2\ell x\in
B(0,ML-2\varepsilon_{\fontsize{8.36pt}{11pt}{\ref{DIDE}}}L+4\ell)$.
\end{itemize}
\end{lemma}

\begin{pf}
First noting that $\delta=\delta_0+4d\varepsilon$, $\delta
_0<2^{-d}S_0$, for
%
%e49 #&#
\begin{equation}
\label{PIM1} \varepsilon_{\fontsize{8.36pt}{11pt}{(\ref{PIM1})}}=\bigl(2^{-d}S_0-
\delta_0\bigr)/(8d)
\end{equation}
$\delta=\delta_0+4d\varepsilon<2^{-d}S_0$ for $\varepsilon\le\varepsilon
_{\fontsize{8.36pt}{11pt}{(\ref{PIM1})}}$. Moreover, for all $x$ such that $2\ell x\in
B(0,ML+4\ell)$,
%
%e50 #&#
\begin{eqnarray}
\label{LBG} \bar G(x,\bar\xi_0)&=&\frac{\mbox{\# of 0's in }\mathcal
{N}(x,K)}{\llvert \mathcal{N}(x,K)\rrvert }\nonumber
\\
&\le&1- \biggl(\frac{\kappa-\varepsilon_{\fontsize{8.36pt}{11pt}{\ref{DIDE}}}-10\varepsilon
}{2\kappa} \biggr)^d S_0
\\
&\le& 1-2^{-d}S_0(1-d\varepsilon_{\fontsize{8.36pt}{11pt}{\ref{DIDE}}}-10d\varepsilon).\nonumber
\end{eqnarray}
Let
%
%e51 #&#
\begin{equation}
\label{PIM2} \varepsilon_{\fontsize{8.36pt}{11pt}{(\ref{PIM2})}}
=\frac{S_0(1-d\varepsilon_{\fontsize{8.36pt}{11pt}{\ref{DIDE}}})-2^d\delta_0}{20\,dS_0+2^{d+2}\,d}>0.
\end{equation}
It is easy to see that $\bar{G}(x,\bar\xi_0)<1-\delta_0-4d\varepsilon
=1-\delta$ for all $\varepsilon\le\varepsilon_{\fontsize{8.36pt}{11pt}{(\ref{PIM2})}}$, which
implies $\bar\omega(\bar G(x,\bar\xi_0))\equiv\omega_0$ for all
$2\ell x\in B(0,ML+4\ell)$. Thus,
\[
\mu_2(x,\bar\xi_0)=\bar\omega\bigl(\bar G(x,\bar
\xi_0)\bigr)f_1(x,\bar\xi_0)-\nu
f_2(x,\bar\xi_0)\ge\omega_0 S_0-
\nu T_0>4\varepsilon_{1} %
\]
for all $2\ell x\in B(0,ML-2\varepsilon_{\fontsize{8.36pt}{11pt}{\ref{DIDE}}}L+4\ell)$. Then for
the infinitesimal mean of type 1:
\begin{eqnarray*}
\mu_1(x,\bar\xi_0)&=&(2L+1)^{-d}\llvert\hat
B_0\rrvert\biggl(\sum_{y\dvtx \hat
B_y\subset\mathcal{N}(x,L)}
f_2(y,\bar\xi_0) \biggr)f_0(x,\bar\xi
_0)\beta
\\
&&{} -[\omega_0+\mu]f_1(x,\bar
\xi_0).
\end{eqnarray*}
Note that
\begin{eqnarray*}
&&\inf_{2\ell x\in B(0,ML+4\ell)}(2L+1)^{-d}\llvert\hat B_0
\rrvert\biggl(\sum_{y\dvtx
\hat B_y\subset\mathcal{N}(x,L)} f_2(y,\bar
\xi_0) \biggr)f_0(x,\bar\xi_0)\beta
\\
&&\qquad \ge2^{-d}(1-3\varepsilon_{\fontsize{8.36pt}{11pt}{\ref{DIDE}}}-10\varepsilon)^dT_0(1-S_0-T_0)
\beta
\\
&&\qquad \ge2^{-d}(1-3d\varepsilon_{\fontsize{8.36pt}{11pt}{\ref{DIDE}}}-10\,d\varepsilon
)T_0(1-S_0-T_0)\beta.
\end{eqnarray*}
Recalling (\ref{PDS}), (\ref{EIDE}) and the definition of $\varepsilon
_1$, let
%
%e52 #&#
\begin{equation}
\label{PIM3} \varepsilon_{\fontsize{8.36pt}{11pt}{(\ref{PIM3})}}=\frac{2^{d-1}\varepsilon
_1}{5\,dT_0(1-S_0-T_0)\beta}.
\end{equation}
For all $\varepsilon\le\varepsilon_{\fontsize{8.36pt}{11pt}{(\ref{PIM3})}}$ and any $2\ell x\in
B(0,ML+4\ell)$,
\[
\mu_2(x,\bar\xi_0)\ge2^{-d}(1-3\,d
\varepsilon_{\fontsize{8.36pt}{11pt}{\ref{DIDE}}})T_0(1-S_0-T_0)\beta-(
\omega_0+\mu)S_0-2\varepsilon_1\ge2
\varepsilon_1. %
\]
Overall, let
\[
\varepsilon_{\fontsize{8.36pt}{11pt}{\ref{PIM}}}=\min\{\varepsilon_{\fontsize{8.36pt}{11pt}{(\ref{PIM1})}},\varepsilon
_{\fontsize{8.36pt}{11pt}{(\ref{PIM2})}},\varepsilon_{\fontsize{8.36pt}{11pt}{(\ref{PIM3})}}\}. %
\]
It satisfies the condition of this lemma by definition.
\end{pf}

Moreover, since that the inequalities for $\bar G$'s in the proof above
are strict and that all other terms in the infinitesimal mean are
continuous, we have:

%
%
%le8.3 #&#
\begin{lemma}\label{Pur}
There is some $\delta_{\fontsize{8.36pt}{11pt}{\ref{Pur}}}>0$ so that for any configuration
$\bar\xi_0'$ with
\[
\bigl\llvert f_i(x,\bar\xi_0)-f_i\bigl(x,
\bar\xi_0'\bigr)\bigr\rrvert\le\delta_{\fontsize{8.36pt}{11pt}{\ref{Pur}}},\qquad
i=1,2; 2\ell x\in B(0,2ML) %
\]
we have:
\begin{itemize}
\item$\mu_1(x,\bar\xi_0')\ge\varepsilon_{1}$ for all $2\ell x\in
B(0,ML+4\ell)$,
\item$\mu_2(x,\bar\xi_0')\ge\varepsilon_{1}$ for all $2\ell x\in
B(0,ML-2\varepsilon_{\fontsize{8.36pt}{11pt}{\ref{DIDE}}}L+4\ell)$
\end{itemize}
for all $\ell=\varepsilon L\le\varepsilon_{\fontsize{8.36pt}{11pt}{\ref{PIM}}}L$.
\end{lemma}

\begin{pf}
First note that for any $x$ such that $2\ell x\in B(0,ML+4\ell)$,
\[
\bigl\llvert\bar G(x,\bar\xi_0)-\bar G\bigl(x,\bar
\xi'_0\bigr)\bigr\rrvert\le2\delta_{\fontsize{8.36pt}{11pt}{\ref{Pur}}}.
\]
Let
%
%e53 #&#
\begin{equation}
\label{Pur1} \delta_{\fontsize{8.36pt}{11pt}{(\ref{Pur1})}}= \bigl[2^{-d}S_0(1-d
\varepsilon_{\fontsize{8.36pt}{11pt}{\ref{DIDE}}}-10\,d\varepsilon_{\fontsize{8.36pt}{11pt}{\ref{PIM}}})-\delta_0-4\,d
\varepsilon_{\fontsize{8.36pt}{11pt}{\ref{PIM}}} \bigr]/4>0.
\end{equation}
For any $\delta_{\fontsize{8.36pt}{11pt}{\ref{Pur}}}\le\delta_{\fontsize{8.36pt}{11pt}{(\ref{Pur1})}}$, recalling
(\ref{LBG}) and (\ref{PIM2}), we have
\begin{eqnarray*}
\bar G\bigl(x,\bar\xi'_0\bigr)&\le&\bar G(x,\bar
\xi_0)+2\delta_{\fontsize{8.36pt}{11pt}{\ref{Pur}}}
\\
&\le&1-2^{-d}S_0(1-d\varepsilon_{\fontsize{8.36pt}{11pt}{\ref{DIDE}}}-10\,d\varepsilon
_{\fontsize{8.36pt}{11pt}{\ref{PIM}}})+2\delta_{\fontsize{8.36pt}{11pt}{(\ref{Pur1})}}
<1-\delta_0-4\,d\varepsilon_{\fontsize{8.36pt}{11pt}{\ref{PIM}}}
\\
&\le& 1-\delta_0-4\,d
\varepsilon=1-\delta
\end{eqnarray*}
which implies that
%
%e54 #&#
\begin{equation}
\label{CD1} \bar\omega\bigl(G\bigl(x,\bar\xi'_0
\bigr)\bigr)\equiv\omega_0
\end{equation}
for all $x$ such that $2\ell x\in B(0,ML+4\ell)$. Furthermore, under
(\ref{CD1}), (\ref{InfM12}) implies that
%
%e55 #&#
%e56 #&#
\begin{eqnarray}
\label{CD2} \qquad\bigl\llvert\mu_2(x,\bar\xi_0)-
\mu_2\bigl(x,\bar\xi_0'\bigr)\bigr\rrvert
&\le&\omega_0\delta_{\fontsize{8.36pt}{11pt}{\ref{Pur}}}+\nu\delta_{\fontsize{8.36pt}{11pt}{\ref{Pur}}}=(
\omega_0+\nu)\delta_{\fontsize{8.36pt}{11pt}{\ref{Pur}}},
\\
 \bigl\llvert\mu_1(x,\bar\xi_0)-
\mu_1\bigl(x,\bar\xi_0'\bigr)\bigr\rrvert&
\le& (2L+1)^{-d}\llvert\hat B_0\rrvert\biggl(\sum
_{y\dvtx \hat B_y\subset\mathcal{N}(x,L)} f_2(y,\bar\xi_0) \biggr
)2\beta \delta_{\fontsize{8.36pt}{11pt}{\ref{Pur}}}\nonumber
\\
\label{CD3} &&{}+(2L+1)^{-d}\llvert\hat B_0\rrvert\biggl(\sum
_{y\dvtx \hat B_y\subset\mathcal
{N}(x,L)} \delta_{\fontsize{8.36pt}{11pt}{\ref{Pur}}} \biggr)f_0(x,\bar
\xi_0)\beta
\\
&&{}+2\delta_{\fontsize{8.36pt}{11pt}{\ref{Pur}}}^2+(\omega_0+\mu)
\delta_{\fontsize{8.36pt}{11pt}{\ref{Pur}}}\nonumber
\end{eqnarray}
for all $x$ such that $2\ell x\in B(0,ML+4\ell)$. Noting that $f_i\le
1$, $\delta_{\fontsize{8.36pt}{11pt}{\ref{Pur}}}\le1$, (\ref{CD3}) can be simplified as
\[
\bigl\llvert\mu_1(x,\bar\xi_0)-\mu_1
\bigl(x,\bar\xi_0'\bigr)\bigr\rrvert\le(2+3\beta+\omega
_0+\mu)\delta_{\fontsize{8.36pt}{11pt}{\ref{Pur}}}. %
\]
Thus, let
\[
\delta_{\fontsize{8.36pt}{11pt}{\ref{Pur}}}=\min\biggl\{\delta_{\fontsize{8.36pt}{11pt}{(\ref{Pur1})}}, \frac
{\varepsilon_1}{2(\omega_0+\nu)},
\frac{\varepsilon_1}{2(2+3\beta
+\omega_0+\mu)} \biggr\}. %
\]
Equations (\ref{CD1})--(\ref{CD3}) show that $\delta_{\fontsize{8.36pt}{11pt}{\ref{Pur}}}$ satisfies
the conditions in our lemma.\vadjust{\goodbreak}
\end{pf}

Both $S_{\mathrm{test}}$ and $T_{\mathrm{test}}$ are Lipchitz with constants
$S_0\varepsilon
_{\fontsize{8.36pt}{11pt}{\ref{DIDE}}}^{-1}$ and $T_0\varepsilon_{\fontsize{8.36pt}{11pt}{\ref{DIDE}}}^{-1}$. Let
$C_{\mathrm{lip}}$ be the max of these two constants. At this point, we are
ready to specify the size of our small box and have the lemma as follows.

%
%le8.4 #&#
\begin{lemma}\label{ExRe}
For $\varepsilon=\varepsilon_{\fontsize{8.36pt}{11pt}{\ref{ExRe}}}$, where
%
%e57 #&#
\begin{equation}
\label{BoxSize} \varepsilon_{\fontsize{8.36pt}{11pt}{\ref{ExRe}}}=\min\biggl\{\frac{\delta_{\fontsize{8.36pt}{11pt}{\ref{Pur}}}\varepsilon_{1}}{16(\beta+\omega_0+\mu)C_{\mathrm{lip}}},
\frac{\varepsilon
_{\fontsize{8.36pt}{11pt}{\ref{PIM}}}}{2} \biggr\},
\end{equation}
$\bar\xi_t$ be the truncated process starting from $\bar\xi_0$. At time
\[
t_{\fontsize{8.36pt}{11pt}{\ref{ExRe}}}=\frac{\delta_{\fontsize{8.36pt}{11pt}{\ref{Pur}}}}{2(\beta+\omega_0+\mu)} %
\]
there is some $C_{\fontsize{8.36pt}{11pt}{\ref{ExRe}}}<\infty$ such that the probability that:
\begin{itemize}
\item$f_1(x,\bar\xi_{t_{\fontsize{8.36pt}{11pt}{\ref{ExRe}}}})\ge f_1(x,\bar\xi_0)+c_{\fontsize{8.36pt}{11pt}{\ref{ExRe}}}$, when $2\ell x\in B(0,ML+4\ell)$,
\item$f_2(x,\bar\xi_{t_{\fontsize{8.36pt}{11pt}{\ref{ExRe}}}})\ge f_2(x,\bar\xi_0)+c_{\fontsize{8.36pt}{11pt}{\ref{ExRe}}}$, when $2\ell x\in B(0,ML-2\varepsilon_{\fontsize{8.36pt}{11pt}{\ref{DIDE}}}L+4\ell)$
\end{itemize}
is greater than $1-C_{\fontsize{8.36pt}{11pt}{\ref{ExRe}}}L^{-d}$ when $L$ is large, where
\[
c_{\fontsize{8.36pt}{11pt}{\ref{ExRe}}}=\frac{\delta_{\fontsize{8.36pt}{11pt}{\ref{Pur}}}\varepsilon_{1}}{8(\beta
+\omega_0+\mu)}. %
\]
\end{lemma}

\begin{pf}
Consider the stopping time
\[
\bar\tau=\min\bigl\{t\dvtx \exists x\dvtx 2\ell x\in B(0,2ML), i=1\mbox{ or }
2, \bigl\llvert
f_i(x,\bar\xi_0)-f_i(x,\bar
\xi_t)\bigr\rrvert> \delta_{\fontsize{8.36pt}{11pt}{\ref{Pur}}}\bigr\}. %
\]
Note that each site in our system flip at a rate no larger than $\beta
+\omega_0+\mu$. According to standard large deviations result as we
used in Lemma \ref{le5.1}, there is some $c_{\fontsize{8.36pt}{11pt}{\ref{MoS}}}, C_{\fontsize{8.36pt}{11pt}{\ref{MoS}}}\in
(0,\infty)$ independent to $L$ such that
\[
P(\bar\tau\le t_{\fontsize{8.36pt}{11pt}{\ref{ExRe}}})\le C_{\fontsize{8.36pt}{11pt}{\ref{MoS}}}\exp\bigl(-c_{\fontsize{8.36pt}{11pt}{\ref{MoS}}}L^d
\bigr)<C_{\fontsize{8.36pt}{11pt}{\ref{MoS}}}L^{-d} %
\]
when $L$ is large. Now consider, $\sigma^2_i(x,\bar\xi_t)$, the
infinitesimal variances of the local densities. According to exactly
the same calculation as we did in Lemma \ref{le3.3}, there is a $C_{\fontsize{8.36pt}{11pt}{\ref{M2tbd}}}<\infty$ such that
%
%e58 #&#
\begin{equation}
\label{InfV12} \sigma^2_i(x,\bar\xi)\le
C_{\fontsize{8.36pt}{11pt}{\ref{M2tbd}}}L^{-d}
\end{equation}
for all $x\in\ZZ^d$, $i=1,2$ and all configurations $\bar\xi$.
Thus, we can again define Dynkin's martingale:
%
%e59 #&#
\begin{equation}
\label{DyM2} \bar M_i(x,t)=f_i(x,\bar
\xi_t)-f_i(x,\bar\xi_0)-\int
_0^t\mu_i(x,\bar
\xi_t)\,dt
\end{equation}
and Lemma \ref{M2tbd} implies that there is a $C_{\fontsize{8.36pt}{11pt}{\ref{M2tbd}}}$ so that
\[
P\Bigl(\sup_{t\le t_{\fontsize{8.36pt}{11pt}{\ref{ExRe}}}}\bigl\llvert\bar{M}_i(x,t)
\bigr\rrvert>c_{\fontsize{8.36pt}{11pt}{\ref{ExRe}}}\Bigr)<C_{\fontsize{8.36pt}{11pt}{\ref{M2tbd}}}L^{-d}.
\]
Consider the event
%
%e60 #&#
\begin{equation}
\label{Aix} A_i(x)=\{\tau>t_{\fontsize{8.36pt}{11pt}{\ref{ExRe}}}\}\cap\Bigl\{\sup
_{t\le t_{\fontsize{8.36pt}{11pt}{\ref{ExRe}}}}\bigl\llvert\bar{M}_i(x,t)\bigr\rrvert
<c_{\fontsize{8.36pt}{11pt}{\ref{ExRe}}}\Bigr\}.
\end{equation}
By definition, there is some $U_{\fontsize{8.36pt}{11pt}{\ref{ExRe}}}<\infty$, independent to
$L$, such that
\[
P\bigl(A_i(x)\bigr)>1-U_{\fontsize{8.36pt}{11pt}{\ref{ExRe}}}L^{-d}\qquad \forall x
\mbox{ s.t. }2\ell x\in B(0,ML+4\ell). %
\]
For any $x$ such that $2\ell x\in B(0,ML+4\ell)$, when $A_1(x)$ holds,
Lemma \ref{Pur} implies that for any $x$
%
%e61 #&#
\begin{equation}
\label{Incr1} f_1(x,\bar\xi_{t_{\fontsize{8.36pt}{11pt}{\ref{ExRe}}}})\ge\int
_0^{t_{\fontsize{8.36pt}{11pt}{\ref{ExRe}}}}\varepsilon_{1}\,dt-c_{\fontsize{8.36pt}{11pt}{\ref{ExRe}}}>c_{\fontsize{8.36pt}{11pt}{\ref{ExRe}}}.
\end{equation}
Similarly, for any $x$ such that $2\ell x\in B(0,ML-2\varepsilon_{\fontsize{8.36pt}{11pt}{\ref{DIDE}}}L+4\ell)$, when $A_2(x)$ holds
%
%e62 #&#
\begin{equation}
\label{Incr162} f_2(x,\bar\xi_{t_{\fontsize{8.36pt}{11pt}{\ref{ExRe}}}})\ge\int
_0^{t_{\fontsize{8.36pt}{11pt}{\ref{ExRe}}}}\varepsilon_{1}\,dt-c_{\fontsize{8.36pt}{11pt}{\ref{ExRe}}}>c_{\fontsize{8.36pt}{11pt}{\ref{ExRe}}}.
\end{equation}
So let
%
%e63 #&#
\begin{equation}
\label{Bevent2} A= \biggl(\bigcap_{x\dvtx 2\ell x\in B(0,ML+4\ell)}A_1(x)
\biggr)\cap\biggl(\bigcap_{x\dvtx2\ell x\in B(0,ML-2\varepsilon_{\fontsize{8.36pt}{11pt}{\ref{DIDE}}}L+4\ell
)}A_2(x)
\biggr).
\end{equation}
The conditions in our lemma are satisfied on the event $A$. Noting that
\[
P(A)\ge1-\sum_{x\dvtx 2\ell x\in B(0,ML+4\ell)}\bigl[ P\bigl(A_1^c(x)
\bigr)+P\bigl(A_2^c(x)\bigr)\bigr]\ge1-
\frac{4M^d}{\varepsilon_{\fontsize{8.36pt}{11pt}{\ref{ExRe}}}^d}U_{\fontsize{8.36pt}{11pt}{\ref{ExRe}}}L^{-d} %
\]
let $C_{\fontsize{8.36pt}{11pt}{\ref{ExRe}}}=4M^dU_{\fontsize{8.36pt}{11pt}{\ref{ExRe}}}/\varepsilon_{\fontsize{8.36pt}{11pt}{\ref{ExRe}}}^d$ and
the proof is complete.
\end{pf}

For any $x\in\ZZ^d$ and any $\xi\in\{0,1,2\}^{\ZZ^d}$, define
$\operatorname{shift}(\xi,x)$ to be the configuration that for any $y\in\ZZ^d$:
\[
\operatorname{shift}(\xi,x) (y)=\xi(y-x). %
\]
Recalling the definition of $C_{\mathrm{lip}}$, on the event $A$, for any
$i=1,\ldots, d$,
\[
\bar\xi_{t_{\fontsize{8.36pt}{11pt}{\ref{ExRe}}}}\ge\operatorname{shift}(\bar\xi_0,\pm2\ell
e_i). %
\]

Monotonicity enables us to restart the construction above from anyone
among the shifts. Note that the success probability of such a
construction is of $1-O(L^{-d})$. So when $L$ is large, with high
probability we can do it for $2d\log L$ times without a failure. This
will give us a ``copy'' of $\bar\xi_0$ at $\pm2\ell(\log L)e_i$ for
each $i$ and will take time $T=(\log L)t_{\fontsize{8.36pt}{11pt}{\ref{ExRe}}}$.

Thus, we can have out block construction as follows: let
\[
\Gamma_x=2\ell(\log L)x+[-\ell\log L,\ell\log L]^d\qquad
\forall x\in\ZZ^d %
\]
and $T_n=nT, n\ge0$. We say $(x,n)$ is wet if
\[
\bar\xi_{T_n}\ge\operatorname{shift}\bigl(\bar\xi_0,2\ell(\log L)x
\bigr). %
\]
From the construction above, we immediately have that $(x,n)$ is wet
then with high probability $(x\pm e_i, n+1)$ are all wet for
$i=1,\ldots, d$.

To check that with high-probability, the block events are
finite-dependent, note that $\bar\xi_t$ is dominated by a branching
random walk with birth rate $\beta$ and initial configuration $\bar
\xi_0$. Lemma \ref{BRW} shows that for any $m>0$
\[
P\bigl(M_k(T_1)\ge(2\beta+m)LT_1\bigr)\le
e^{-mT_1}\bigl\llvert B\bigl(0,(M+1)L\bigr)\bigr\rrvert%
\]
when $L$ is large enough, where $M_k(t)$ is the largest $k$th
coordinate among the occupied sites at time $t$. Noting that $T_1=(\log
L)t_{\fontsize{8.36pt}{11pt}{\ref{ExRe}}}$ and that the choice of $t_{\fontsize{8.36pt}{11pt}{\ref{ExRe}}}$ is
independent to the choice of $L$, let
%
%e64 #&#
\begin{equation}
\label{DSize} m=(d+1)/t_{\fontsize{8.36pt}{11pt}{\ref{ExRe}}}.
\end{equation}
We can control the probability that $\bar\xi_t$ wanders too far as follows:
\[
P_{\bar\xi_0} \Bigl(\max_{t\in[0,T_0]}\bigl\{\llVert x\rrVert\dvtx
\bar\xi_t(x)\neq0\bigr\}\ge(2\beta+m)LT_1 \Bigr)\le2\,d
e^{-mT_1}\bigl\llvert B\bigl(0,(M+1)L\bigr)\bigr\rrvert. %
\]
Noting that
\[
e^{-mT_1}=e^{-(\log L)t_{\fontsize{8.36pt}{11pt}{\ref{ExRe}}}(d+1)/t_{\fontsize{8.36pt}{11pt}{\ref{ExRe}}}}=e^{-(d+1)\log L}=L^{-d-1} %
\]
we have
%
%e65 #&#
\begin{equation}
\label{FDep} P_{\bar\xi_0} \Bigl(\max_{t\in[0,T_0]}\bigl\{
\llVert x\rrVert\dvtx \bar\xi_t(x)\neq0\bigr\}\ge(2\beta+m)LT_1
\Bigr)=L^{-d-1} O\bigl(L^d\bigr)\rightarrow0
\end{equation}
as $L\rightarrow\infty$. Thus, noting that $\ell=\varepsilon_{\fontsize{8.36pt}{11pt}{\ref{ExRe}}}L$, let
\[
R=\frac{(2\beta+m)LT_1}{2\ell\log L}=\frac{(2\beta+m)t_{\fontsize{8.36pt}{11pt}{\ref{ExRe}}}}{2\varepsilon_{\fontsize{8.36pt}{11pt}{\ref{ExRe}}}} %
\]
which is a finite constant independent to the choice of $L$. As $L$
goes large, we have that with high probability $\bar\xi_t$ cannot
exit the following finite union of blocks by time $T_1$:
\[
\bigcup_{x\dvtx \llVert x\rrVert \le R}\Gamma_x %
\]
which implies that the block events we constructed has finite range of
dependence. Then again according to standard block argument in \cite
{DStF} and \cite{HT}, we complete the proof of survival for $\bar\xi
_t$ and this implies Theorem \ref{SurHD}.

%\begin{appendix}
%\section{}
%\end{appendix}

% zodis "Acknowledgments" paliekamas pagal autoriu
\section*{Acknowledgment} The authors would like to thank two anonymous
referees for many comments that helped to improve the article.

%\begin{supplement}[id=suppA]
%\sname{Supplement A}
%\stitle{}
%\slink[doi]{10.1214/00-AAPXXXXSUPP} %[doi,text={...}] - jei reikia
%suskaldyti doi
%\sdatatype{.pdf}
%\sfilename{aapXXXX\_supp.pdf}
%\sdescription{}
%\end{supplement}

% imsref loaded by linak, 2014-12-22 16:04:01
% imsref loaded by linak, 2014-12-29 12:36:47
% imsref loaded by linak, 2015-01-22 12:50:20

\printaddresses
\end{document}